\documentclass[12pt]{amsart} 
\usepackage[mathscr]{eucal}
\usepackage{amsmath,amsfonts}

\parskip=\smallskipamount

\newtheorem{theorem}{Theorem}[section]

\newtheorem{proposition}[theorem]{Proposition}

\newtheorem{lemma}[theorem]{Lemma}
\newtheorem{example}[theorem]{Example}

\newtheorem{remark}[theorem]{Remark}

\makeatletter
\@addtoreset{equation}{section}
\makeatother

\newcommand{\CC}{{\mathbb C}}
\newcommand{\NN}{{\mathbb N}}

\newcommand{\cA}{{\mathcal A}}
\newcommand{\cB}{{\mathcal B}}

\newcommand{\cD}{{\mathcal D}}

\newcommand{\cF}{{\mathcal F}}
\newcommand{\cG}{{\mathcal G}}
\newcommand{\cH}{{\mathcal H}}
\newcommand{\cK}{{\mathcal K}}
\newcommand{\cL}{{\mathcal L}}

\newcommand{\cN}{{\mathcal N}}
\newcommand{\cP}{{\mathcal P}}

\newcommand{\cU}{{\mathcal U}}

\newcommand{\fI}{{\mathfrak I}}

\newcommand{\lhk}{\cL(\cH,\cK)}

\newcommand{\lh}{\cL(\cH)}
\newcommand{\lk}{\cL(\cK)}

\newcommand{\rank}{\hbox{\rm{rank}}\,}
\newcommand{\sgn}{\hbox{\rm{sgn}}\,}
\newcommand{\Exp}{\hbox{\rm{Exp}}\,}
\newdimen\expt
\def\boxit#1{\setbox0\hbox{$\displaystyle{#1}$}
      \hbox{\lower.4\expt
 \hbox{\lower3\expt\hbox{\lower\dp0
      \hbox{\vbox{\hrule height.4\expt
 \hbox{\vrule width.4\expt\hskip3\expt
      \vbox{\vskip3\expt\box0\vskip2\expt}%
 \hskip3\expt\vrule width.4\expt}\hrule height.4\expt}}}}}}
\begin{document}
\pagestyle{plain}

\bigskip

\title 
{Representations of Hermitian Kernels \\
by Means of Kre\u{\i}n spaces \\
II. Invariant Kernels} 
\author{T. Constantinescu} \author{A. Gheondea} 

\address{Department of Mathematics \\
  University of Texas at Dallas \\
  Box 830688, Richardson, TX 75083-0688, U.S.A.}
\email{\tt tiberiu@utdallas.edu} 
\address{Institutul de Matematic\u{a}
  al Academiei Rom\^{a}ne \\
  C.P. 1-764, 70700 Bucure\c{s}ti, Rom\^{a}nia} 
\email{\tt gheondea@imar.ro}

\begin{abstract}
In this paper we study hermitian kernels invariant under the action of a 
semigroup with involution. We characterize those hermitian kernels 
that realize 
the given action by bounded operators on a Kre\u{\i}n space.
Applications to the GNS representation of $*$-algebras associated to hermitian
functionals are given.
We explain the key role played by the Kolmogorov decomposition in
the construction of Weyl exponentials associated to an indefinite inner
product and in the dilation theory of hermitian maps on $C^*$-algebras. 
\end{abstract}

\maketitle

\section{Introduction}

The Hilbert space $\cH$ associated to a positive definite kernel $K$ 
is an abstract version of the $L^2$ space associated to a positive 
measure and the Kolmogorov decomposition of $K$ gives a useful 
expansion of the elements of $\cH$ in terms of a geometrical
model of a stochastic process with covariance kernel $K$. Therefore, it
is quite natural to seek similar constructions for an 
arbitrary kernel. While the decomposition into a real and
an imaginary part can be realized without difficulties, the 
study of hermitian kernels is no longer straightforward.
This was shown in the work of L.~Schwartz \cite{Sch},
where a characterization of the hermitian kernels admitting
a Jordan decomposition was obtained in terms of a boundedness
condition that we call the Schwartz condition (the statement (1) 
of Theorem~\ref{kolmo} in Section~2 below). A key difficulty of the theory
was identified in \cite{Sch} in the lack of uniqueness of the 
associated reproducing kernel spaces.

It was shown in \cite{CG1} that the Schwartz condition is also
equivalent to the existence of a Kolmogorov decomposition, 
while the uniqueness of the Kolmogorov decomposition was
characterized in spectral terms (Theorem~\ref{kolmo} and, 
respectively, Theorem~\ref{unic} in Section~2). 
The purpose of this paper 
is to continue these investigations by considering
hermitian kernels with additional symmetries given by the 
action of a semigroup. The main result gives a characterization
of those hermitian kernels that produce a representation
of the action by bounded operators on a certain Kre\u{\i}n space.
It turns out that such a result has many applications and in this paper
we discuss GNS representations on inner product spaces and dilation
theory.

The paper is organized as follows. In Section~2 we review the concept 
of induced Kre\u{\i}n space and we show its key role in the 
construction of Kolmogorov decompositions as described in \cite{CG1}.
A new result is added here in connection with a lifting
property for induced Kre\u{\i}n spaces that is related to an important
inequality of M.~G.~Kre\u{\i}n. 
Theorem~\ref{contraexemplu} gives an example of an induced
Kre\u{\i}n space without the lifting property, adding one more
pathology to the study of hermitian kernels. Incidently, this
result answers negatively a question raised in \cite{DR}.
In Section~3 we prove the main result of the paper.
We consider the action of a semigroup on a hermitian kernel
and Theorem~\ref{baza} gives the conditions that insure
the representation of this action as a semigroup of bounded 
operators on a Kre\u{\i}n space. 
We also address the 
uniqueness property of such representations. While the case
of the trivial semigroup with one element is settled in \cite{CG1}
(Theorem~\ref{unic} in Section~2) and Theorem~\ref{pontr} 
gives another partial answer, the general
case remains open. The proof used for the 
trivial semigroup cannot be easily extended precisely because
Theorem~\ref{contraexemplu} is true.
In Section~4 we analyze the case when the projective representation 
given by Theorem~\ref{baza} is fundamentally reducible or,
equivalently, it is similar with a projective Hilbert space
representation, a question closely related to other similarity problems
and uniformly bounded representations. 
Section~5 shows as an example how to use
the Kolmogorov decomposition in order to obtain Weyl exponentials
associated to an indefinite inner product.
In Section~6 we apply our results to questions 
concerning GNS representations of $*$-algebras on 
Kre\u{\i}n spaces. The whole issue is motivated by the lack
of positivity in some models in local quantum field theories.
We relate these questions to properties of Kolmogorov 
decompositions so that we can characterize the existence
(Theorem~\ref{starex}), the uniqueness (Theorem~\ref{starun}),
and the boundedness of GNS data (Theorem~\ref{starmar}).
 In Section~7 we obtain 
a version of the Stinespring representation 
theorem for hermitian maps. It turns out that this result
not only extends the classical Stinespring theorem but it is also
related to the Wittstock representation of completely bounded
maps. The main point of this discussion is to show that
the Schwartz condition and its corresponding version in 
Theorem~\ref{baza} are related to the concept 
of completely bounded map.

\section{Preliminaries}  

We briefly review the concept of a Kolmogorov decomposition
for hermitian kernels.
The natural framework to deal with these kernels is that
of Kre\u{\i}n spaces. We recall first some definitions 
and a few items of notation.
An indefinite inner product space
$(\cH,[\cdot,\cdot])$ is called {\it Kre\u{\i}n space}
provided that there exists 
a positive inner product $\langle\cdot,\cdot\rangle $
on $\cH $ turning $(\cH, \langle\cdot,\cdot\rangle )$ 
into a Hilbert space
such that $[\xi ,\eta ]=\langle J\xi ,\eta \rangle $, 
$\xi ,\eta \in \cH $,
for some symmetry $J$ ($J^*=J^{-1}=J$) on $\cH$. Such a 
symmetry $J$ is called a {\it fundamental symmetry}. The norm 
$\|\xi\|^2=\langle \xi,\xi\rangle$ is called a {\it unitary norm}.
The underlying Hilbert space topology of $\cK$ is called the {\em
strong topology} and does not depend on the choice of the fundamental symmetry.

For two Kre\u{\i}n spaces $\cH$ and $\cK$ 
we denote by $\lhk$ the set 
of linear bounded operators from $\cH $ to $\cK$.
For $T\in \lhk$ we denote by $T^{\sharp }$ the adjoint 
of $T$ with respect to $[\cdot,\cdot]$. 
We say that $A\in \lh$ is a {\it selfadjoint 
operator} if $A^{\sharp }=A$. A possibly unbounded operator $V$ between 
two Kre\u\i n spaces is called {\it isometric} if $[V\xi,V\eta]=[\xi,\eta]$
for all $\xi,\eta$ in the domain of $V$.
Also, we say that the operator $U\in \lh$ is
{\it unitary} if $UU^{\sharp }=U^{\sharp }U=I$, where 
$I$ denotes the identity operator on $\cH$.
The notation 
$T^*$ is used for the adjoint of $T$ with respect to the
positive inner product $\langle\cdot,\cdot\rangle $.

\subsection{Kre\u\i n Spaces Induced by Selfadjoint Operators}

Many difficulties in dealing with operators on 
Kre\u{\i}n spaces are caused by the 
lack of a well-behaved factorization 
theory. The concept of induced space
turned out to be quite useful in this direction.
Thus, for a selfadjoint operator $A$ in $\lh$
we define a new inner product $[\cdot ,\cdot ]_A$
on $\cH$ by the formula
\begin{equation}\label{subA}
[\xi ,\eta ]_A=[A\xi ,\eta ], \quad \xi ,\eta \in \cH ,
\end{equation} 
and a pair $(\cK, \Pi )$
consisting of a Kre\u{\i}n space $\cK$
and a bounded operator $\Pi \in \lhk$ 
is called a {\it Kre\u{\i}n space induced}
by $A$ provided that $\Pi $ has dense range
and the relation 
\begin{equation}\label{indus}
[\Pi \xi ,\Pi \eta ]_{\cK }=[\xi ,\eta ]_A
\end{equation} 
holds for all $\xi ,\eta \in \cH$, where 
$[\cdot ,\cdot ]_{\cK }$ denotes the indefinite inner 
product on $\cK $.
One well-known example is obtained in the following way.

\begin{example}\label{hindus}{\em
Let $J$ be a fundamental symmetry on $\cH$
and let $\langle \cdot ,\cdot  \rangle _J$ be the associated positive
inner product turning $\cH$ into a Hilbert space.
Then $JA$ is a selfadjoint operator on this Hilbert space
and let $\cH _-$ and $\cH _+$ be the spectral subspaces of 
$JA$ corresponding to $(-\infty ,0)$ and, respectively,
$(0, \infty )$.
We obtain the decomposition 
$$\cH =\cH _-\oplus \ker A\oplus \cH _+.$$
Note that $(\cH _-, -[\cdot ,\cdot ]_A)$ 
and $(\cH _+, [\cdot ,\cdot ]_A)$
are positive inner product spaces and hence they 
can be completed to the Hilbert spaces
$\cK _-$ and, respectively, $\cK_+$. Let $\cK _A$ be the 
Hilbert direct sum of 
$\cK _-$ and $\cK_+$ and denote by 
$\langle \cdot ,\cdot \rangle _{\cK _A}$
the positive inner product on $\cK _A$.
Define
$$
J_A(k_-\oplus k_+)=-k_-\oplus k_+
$$
for $k_-\in \cK _-$ and $k_+\in \cK _+$.
We can easily check that $J_A$ is a symmetry 
on $\cK _A$ and then the inner product
$$[k,k']_{\cK _A}=
\langle J_Ak,k'\rangle _{\cK _A}$$
turns $\cK _A$ into a Kre\u{\i}n space.
The map $\Pi _A:\cH \rightarrow \cK _A$
is defined by the formula
$$\Pi _A\xi =[P_{\cH _-}\xi ]\oplus [P_{\cH _+}\xi ],$$
where $\xi \in \cH$, $P_{\cH _{\pm}}$
denotes the orthogonal projection 
of the Hilbert space $(\cH ,\langle \cdot ,\cdot \rangle _J)$
onto the subspace $\cH _{\pm }$, and $[P_{\cH _{\pm }}\xi ]$
denotes the class of $P_{\cH _{\pm }}\xi $ in $\cK _{\pm }$.
Then one checks that $(\cK _A,\Pi _A)$
is a Kre\u{\i}n space induced by $A$.
In addition, if $JA=S_{JA}|JA|$
is the polar decomposition of $JA$, then we 
note that 
\begin{equation}\label{JPI}
J_A\Pi _A=\Pi _AS_{JA}.
\end{equation}
}
\qed\end{example}

This example proved to be very useful since it is
accompanied by a good property concerning the 
lifting of operators, as shown by a classical
result of M.G. Kre\u{\i}n, \cite{Kr}.
The result was rediscovered by W.J. Reid
\cite{Rd}, P.D. Lax \cite{Lx}, 
and J. Dieudonn\'e \cite{Dd}. The indefinite version 
presented below was proved in \cite{DLS} by using a $2\times2$ 
matrix construction that reduces the proof to the positive definite case.

\begin{theorem}\label{krein}
Let $A$ 
and $B$
be bounded selfadjoint operators 
on the Kre\u{\i}n
spaces $\cH _1$ and $\cH _2$.
Assume that the operators $T_1\in \cL(\cH _1,\cH _2)$
and $T_2\in \cL(\cH _2,\cH _1)$
satisfy the relation $T_2^{\sharp }A=BT_1$.
Then there exist (unique) operators
$\tilde T_1\in \cL(\cK _A,\cK _B)$ and 
$\tilde T_2\in \cL(\cK _B,\cK _A)$
such that 
$\tilde T_1\Pi _A=\Pi _BT_1$, 
$\tilde T_2\Pi _B=\Pi _AT_2$,
and $[\tilde T_1f,g]_{\cK _B}=[f,\tilde T_2g]_{\cK _A}$
for all $f\in \cK _A$, $g\in \cK _B$.
\end{theorem}

Theorem~\ref{krein} will be used 
in an essential way in the 
proof of the main result of the next section and it is also
 related to the uniqueness property  of a Kolmogorov
decomposition for invariant hermitian kernels. 
For these reasons we discuss one more 
question related to this result, namely whether 
this lifting property holds for other induced Kre\u\i n spaces. 
More precisely, two
Kre\u\i n spaces $(\cK_i,\Pi_i)$, $i=1,2$, induced by the same selfadjoint
operator $A\in\cL(\cH)$ are {\it unitarily equivalent} if there exists
a unitary operator $U$ in $\cL(\cK_1,\cK_2)$ such that $U\Pi_1=\Pi_2$.
Theorem~2.8 in \cite{CG1} shows that there exist selfadjoint
operators with the property that not all of their 
induced Kre\u\i n spaces are unitarily equivalent.

Let $(\cK,\Pi )$ be a Kre\u{\i}n space induced by $A$.
We say that $(\cK,\Pi )$ has the {\it lifting property} 
if for any pair of operators $T,S\in \lh$ satisfying the relation
$AT=SA$ there exist unique operators $\tilde T, \tilde S\in \lk$
such that $\tilde T\Pi =\Pi T$, $\tilde S\Pi =\Pi S$.
From Theorem~\ref{krein} it follows that the induced
Kre\u\i n space $(\cK_A,\Pi_A)$
constructed in Example~\ref{hindus} has the lifting property, as
do all the others which are unitarily equivalent to it.
However, as the following result shows, this is not true for
all induced Kre\u\i n spaces of $A$.

\begin{theorem}\label{contraexemplu}
There exists a selfadjoint operator 
that has an induced Kre\u{\i}n space without the 
lifting property. 
\end{theorem}
\begin{proof}
Let $\cH_0$ be an infinite dimensional Hilbert space and $A_0$ is a
bounded selfadjoint operator in $\cH_0$ such that $0\leq A_0\leq I$,
$\ker A_0=0$, and the spectrum of $A_0$ accumulates to
$0$, equivalently, its range is not closed. 
Consider the Hilbert space $\cH=\cH_0\oplus\cH_0$ as well as the
bounded selfadjoint operator
\begin{equation}\label{dirac}
A=\left[\begin{array}{cc} A_0 & 0 \\ 0 &
-A_0\end{array}\right].
\end{equation}
Let $\cK$ be the Hilbert space $\cH$ with the indefinite inner prodcut
$[\cdot,\cdot]$ defined by the symmetry
$$J=\left[\begin{array}{cc} I & 0 \\ 0 &
-I\end{array}\right].$$

Consider the operator $\Pi\in\cL(\cH,\cK)$
\begin{equation}\label{aul}
\Pi=\left[\begin{array}{cc} I & -(I-A_0)^{1/2} \\ 
(I-A_0)^{1/2}  & - I\end{array}\right].
\end{equation}
It is a straightforward calculation to see that $\Pi^* J\Pi=A$ and, by
performing a Frobenius-Schur factorization, it follows that $\Pi$ has
dense range. Thus, $(\cK,\Pi)$ is a Kre\u{\i}n space induced by $A$
and we show that it does not have the lifting property.

Let $T$ be an operator in $\cL(\cH)$ such that, with respect to its
$2\times 2$ block-matrix representation, all its entries $T_{ij}$,
$i,j=1,2$, commute with $A_0$. Define the operator 
$S=JTJ$ and note that $AT=SA$. 

Let us assume that there exists a bounded operator $\widetilde
T\in\cL(\cK)$ such that $\widetilde T\Pi=\Pi T$. Then, there exists
the constant $C=\|\widetilde T\|_{\cK}<\infty$ such that
$$\|\Pi T\xi\|\leq C\|\Pi \xi\|,\quad \xi\in\cH,$$ 
or, equivalently, that
\begin{equation}\label{cinci}
T^*HT\leq C^2 H,
\end{equation}
where
$$H=\left[\begin{array}{cc} 2-A_0 & -2(I-A_0)^{1/2} \\ 
-2(I-A_0)^{1/2}  & 2-A_0 \end{array}\right].$$

Taking into account that $A_0$ commutes with all the other operator entries 
involved in (\ref{cinci}),
it follows that the inequality (\ref{cinci}) is equivalent to
\begin{equation}\label{sase} 
T^* \left[\begin{array}{cc} I & -\Delta \\
-\Delta & I \end{array}\right]T\leq C^2 
\left[\begin{array}{cc} I & -\Delta \\
-\Delta & I \end{array}\right],
\end{equation}
where we denoted
$$\Delta=2(I-A_0)^{1/2}(2-A_0)^{-1}.$$
Note that, by continuous functional calculus, $\Delta$ is an operator
in $\cH$ such that $0\leq \Delta\leq I$ and its spectrum accumulates
to $1$. 

The use of the Frobenius-Schur
factorization
\begin{equation}\label{fs}\left[\begin{array}{cc} I & -\Delta \\
-\Delta & I \end{array}\right]=\left[\begin{array}{cc} I & 0 \\
-\Delta & I \end{array}\right]\left[\begin{array}{cc} I & 0 \\
0 & I-\Delta^2 \end{array}\right]\left[\begin{array}{cc} I & -\Delta \\
0 & I \end{array}\right]\end{equation}
suggests to take
$$T=\left[\begin{array}{cc} I & -\Delta \\
0 & I \end{array}\right],$$ 
and this choice is consistent 
with our assumption that all its entries commute with $A_0$.

Since $T$ is bounded invertible, from 
(\ref{sase}) we get
$$\left[\begin{array}{cc} I & -\Delta \\
-\Delta & I \end{array}\right]\leq C^2\left[\begin{array}{cc} I & 0 \\
0 & I-\Delta^2 \end{array}\right].$$
Looking at the lower right corners of
the matrices in the previous inequality we get $I\leq C^2(I-\Delta^2)$
which yields a contradiction since the spectrum of the operator
$I-\Delta^2$ accumulates to $0$.
\end{proof}

\begin{remark} {\em
$(1)$ Incidently, the example in Theorem~\ref{contraexemplu} 
can be used to answer the following question raised in \cite{DR}, Lecture~6A:
let $A$ be a selfadjoint operator on a 
Kre\u{\i}n space $\cH$, and construct a factorization 
$A=DD^{\sharp }$, where $D\in \cL (\cD,\cH)$ is an one-to-one operator.
If $X\in \lh$ and $XA$ is selfadjoint, does there 
exist a (unique) selfadjoint operator $Y\in \cL(\cD)$
such that $XD=DY$?

We show that the answer to this question is negative.
Indeed, an operator $D$ as above produces
the induced Kre\u{\i}n space $(\cD,D^\sharp)$ for $A$. 
Let $A$ be the operator defined by \eqref{dirac}.
Let us take
$$T=\left[\begin{array}{cc} I & I \\ -I & I \end{array}\right].$$ 
One checks that
$AT=T^*A$. Define $X=T^*$, then $XA $ is selfadjoint.
If $Y\in \cL(\cD)$ exists such that $XD=DY$, then 
$Y^*D^*=D^*T$ and a similar reasoning as in the proof of
Theorem~\ref{contraexemplu} shows that from
(\ref{sase}) and (\ref{fs}) we get
$2(\Delta ^3+\Delta ^2-\Delta +I)\leq C^2(I-\Delta ^2)$, which is 
impossible since the spectrum of the operator from the left side is
bounded away from $0$.

$(2)$ One might ask whether another additional assumption on the
operator $T$ that is frequently used in applications,
namely that $T$ is $A$-isometric, could inforce the lifting property. 
To see that this is not the case, let us take 
$$T=\frac{2}{\sqrt{3}}
\left[\begin{array}{cc} I & -\frac{1}{2}I \\ \frac{1}{2}I & -I
\end{array}\right].$$
It is easy to prove that $T^*AT=A$, that is,
$T$ is $A$-isometric. Noting that $T$ is boundedly invertible, this
corresponds to $S=T^{*-1}$. As before, from (\ref{sase}) and (\ref{fs}) we get
$\frac{4}{3}(-\Delta ^3+\frac{15}{4}\Delta ^2-3\Delta +\frac{5}{4}I)
\leq C^2(I-\Delta^2)$.
But this is again contradictory since the spectrum of the operator
from the left side is bounded away from $0$.
}
\end{remark}

\subsection{Kolmogorov Decompositions of Hermitian Kernels}

We can use the concept of induced space in order
to describe the Kolmogorov decomposition of a 
hermitian kernel.
Let $X$ be an arbitrary set.
A mapping $K$ defined on $X\times X$ with values in 
$\lh$, where $(\cH,[\cdot ,\cdot ]_{\cH})$ is a 
Kre\u{\i}n space, 
is called a {\it hermitian kernel on} 
$X$ if
$K(x,y)=K(y,x)^{\sharp }$ for all $x,y\in X$.

Let ${\cF}_0(X,\cH)$ denote the vector space of $\cH$-valued functions
on $X$ having finite support.
We associate to $K$ an inner product on 
${\cF}_0(X,\cH)$
by the formula:
\begin{equation}\label{fega}
[f,g]_K=\sum_{x,y\in X}
[K(x,y)f(y),g(x)]_{\cH}
\end{equation}
for $f,g\in{\cF}_0(X,\cH)$.
We say that the hermitian kernel $L:X\times X\rightarrow \lh$
is {\it positive definite} if the inner product
$[\cdot ,\cdot ]_L$ associated to $L$ 
by the formula \eqref{fega} is positive. 
On the set of hermitian kernels on $X$ with values in $\lh$
we also have a natural partial order defined as follows:
if $A$, $B$ are hermitian kernels, then $A\leq B$ means
$[f,f]_A\leq [f,f]_B$ for all $f\in {\cF}_0(X,\cH)$.
Following L. Schwartz \cite{Sch}, we say that two
positive definite kernels $A$ and $B$ are {\it disjoint}
if for any positive definite kernel $P$ such that
$P\leq A$ and $P\leq B$ it follows that $P=0$.
A {\it Kolmogorov decomposition}
of the hermitian kernel $K$ is a pair $(V;\cK)$, 
where $\cK$ is a Kre\u{\i}n space
and $V=\{V(x)\}_{x\in X}$ is a family of bounded operators
$V(x)\in\lhk $ such that
$K(x,y)=V(x)^{\sharp }V(y)$ for all $x,y\in X$,
and 
the closure of $\bigvee _{x\in X}V(x){\cH}$
is $\cK$
(\cite{Ko}, \cite{PS}, \cite{EL}). Note that here and througout this
paper $\vee$ stands for the linear manifold generated by some set,
without taking any closure.

The next result, obtained in \cite{CG1}, settles the question
concerning the existence of a Kolmogorov decomposition  for a given 
hermitian kernel.

\begin{theorem}\label{kolmo}
Let $K:X\times X\rightarrow \lh$
be a hermitian kernel. The following assertions are equivalent:

{\rm (1)} There exists 
a positive definite kernel $L:X\times X\rightarrow \lh$
such that $-L\leq K\leq L$.

{\rm (2)} $K$ has a Kolmogorov decomposition.
\end{theorem}
The condition $-L\leq K\leq L$ of the 
previous result appeared in the work of 
L. Schwartz \cite{Sch} concerning the structure of hermitian
kernels. We will call it the {\it Schwartz condition}. 
It is proved in \cite {Sch} that this condition
is also 
equivalent to the Jordan decomposition of $K$, 
which means that 
the kernel $K$ is a difference of two disjoint 
positive definite kernels.
It is convenient for our purpose to review
the construction of the Kolmogorov decomposition.
We assume that there exists 
a positive definite kernel $L:X\times X\rightarrow \lh$
such that $-L\leq K\leq L$.
Let ${\cH}_L$ be the Hilbert space
obtained by the completion of the quotient space
${\cF}_0(X,\cH)/{\cN}_L$ with respect to $[\cdot ,\cdot ]_L$, 
where ${\cN}_L=\{f\in {\cF}_0(X,\cH)\mid [f,f]_L=0\}$
is the isotropic subspace of the inner product space
$({\cF}_0(X,\cH),[\cdot ,\cdot ]_L)$.
Since $-L\leq K\leq L$ is equivalent to
\begin{equation}\label{schwarz}
|[f,g]_K|\leq [f,f]_L^{1/2}\, [g,g]_L^{1/2}
\end{equation}
for all $f,g\in{\cF}_0(X,\cH)$
(see Proposition~38, \cite{Sch}), 
it follows that ${\cN}_L$ is a subset of the 
isotropic subspace ${\cN}_K$ of the inner product space
$({\cF}_0(X,{\cH}),[\cdot ,\cdot ]_K)$.
Therefore, $[\cdot ,\cdot ]_K$ uniquely induces an inner product on 
${\cH}_L$, still denoted by $[\cdot ,\cdot ]_K$, such that
\eqref{schwarz} holds
for $f,g\in{\cH}_L$.
By the Riesz representation theorem we obtain a 
selfadjoint contractive operator $A_L\in{\cL}({\cH}_L)$,
referred to as the {\it Gram}, or {\it metric operator} of $K$ 
with respect to $L$,  such that 
\begin{equation}\label{gram}
[f,g]_K=[A_Lf,g]_L
\end{equation}
for all $f,g\in {\cH _L}.$
Let $(\cK _{A_L},\Pi _{A_L})$ be the  Kre\u{\i}n space induced by $A_L$
given by Example~\ref{hindus}. For $\xi \in \cH$ and $x\in X$, 
we define the element $\xi _x\in {\cF}_0(X,\cH)$
by the formula:
\begin{equation}\label{xi}
\xi _x(y)=\left\{\begin{array}{ll}
\xi, & y=x \\
0, & y\ne x.
\end{array}\right.
\end{equation}
Then we define
\begin{equation}\label{v}
V(x)\xi =\Pi _{A_L}[\xi _x],
\end{equation}
where $[\xi _x]$ denotes the class of $\xi _x$ in ${\cH}_L$ and 
it can be verified that $(V;\cK _{A_L})$ 
is a Kolmogorov decomposition of 
the kernel $K$.

We finally review 
the uniqueness property of the Kolmogorov decomposition.
Two Kolmogorov decompositions $(V_1,\cK _1)$ and
$(V_2,\cK _2)$ of the same hermitian
kernel $K$ are {\it unitarily
equivalent} if there exists a unitary operator 
$\Phi\in\cL(\cK _1,\cK _2)$
such that for all $x\in X$ we have $V_2(x)=\Phi V_1(x)$.
The following result was obtained in \cite{CG1}.

\begin{theorem}\label{unic}
Let $K$ be a hermitian kernel which has Kolmogorov
decompositions. The following assertions are equivalent:\smallskip

{\rm (1)} All Kolmogorov decompositions of $K$ are unitarily equivalent.

{\rm (2)} For each positive definite
kernel $L$ such that $-L\leq K\leq L$, there exists $\epsilon>0$ such that
either $(0,\epsilon)\subset\rho(A_L)$ or $(-\epsilon,0)\subset\rho(A_L)$,
where $A_L$ is the Gram operator of $K$ with respect to $L$.
\end{theorem}

\section{Invariant hermitian kernels}

In this section we study properties
of the Kolmogorov decompositions of 
hermitian kernels with additional symmetries. Let
$S$ be a unital semigroup and $\phi $ an action
of $S$ on the set $X$, this means that $\phi\colon S\times
X\rightarrow X$, $\phi(a,\phi(b,x))=\phi(ab,x)$ for all $a,b\in
S$, $x\in X$, and $\phi (e,x)=x$, where $e$ denotes the unit
element of $S$. 
We are interested in those kernels $K$ on $X$ 
assumed to satisfy a certain invariance property
with respect to the action $\phi $ because this leads to the 
construction of a representation of $S$ on the space of a Kolmogorov
decomposition of $K$.
This kind of construction is well-known for a positive definite kernel
(it just extends the construction 
of the regular representation, see for instance, \cite{PS}), 
but for the Kre\u{\i}n space 
setting the question concerning the boundedness of the  
representation operators is more delicate. 
It is the goal of this section to deal with this matter
in a more detailed way. 

We now introduce additional notation and terminology.
Let $\alpha $ be a complex-valued function on $S\times X$ such that
$\alpha (a,x)\ne 0$. Assume that $\alpha $ satisfies the
following relation:
\begin{equation}\label{cob}
\alpha (ab,x)\overline{\alpha (ab,y)}=
\alpha (a,\phi (b,x))\overline{\alpha (a,\phi (b,y))}
\alpha (b,x)\overline{\alpha (b,y)}
\end{equation}
for all $x,y\in X$. This implies that 
$$\sigma (a,b)=\alpha (a,\phi (b,x))^{-1}\alpha (b,x)^{-1}
\alpha (ab,x)$$
does not depend on $x$; moreover, $|\sigma (a,b)|=1$,
and $\sigma $ has the $2$-cocycle property:
\begin{equation}\label{cociclu}
\sigma (a,b)\sigma (ab,c)=\sigma (a,bc)\sigma (b,c) 
\end{equation}
for all $a,b,c\in S$ (see \cite{PS}, Lemma~2.2). 

For $a\in S$ define the mapping
\begin{equation}\label{actpsi}
\psi ^0_a(\xi _x)=\alpha (a,x)^{-1}\xi _{\phi (a,x)}(=
(\alpha (a,x)^{-1}\xi) _{\phi (a,x)}),
\end{equation}
where $\xi_x$ is defined as in \eqref{xi}. Note that
each element $h$ of ${\cF}_0(X,\cH)$ can be uniquely
written as a finite sum $h=\sum _{k=1}^n\xi ^k_{x_k}$
for vectors $\xi^1,\ldots,\xi^k\in \cH$ and distinct
elements $x_1$, $x_2$, $\ldots $, $x_n$ in $X$. Then
the map $\psi ^0_a$ can be extended by linearity
to a linear map $\psi _a$, from ${\cF}_0(X,\cH)$
into ${\cF}_0(X,\cH)$,
$$\psi _a(\sum _{k=1}^n\xi ^k_{x_k})=
\sum _{k=1}^n\psi ^0_a(\xi ^k_{x_k}).$$

We say that a positive definite kernel $L$ is
$\phi$-{\it bounded} provided that for all $a\in S$, 
$\psi _a$ is bounded with respect to the seminorm
$[\cdot ,\cdot ]_L^{1/2}$ induced by $L$ on ${\cF}_0(X,\cH)$ and let
us denote by
${\cB}_{\phi }^+(X,\cH)$ the set of positive definite $\phi$-bounded
kernels on $X$ with values in $\lh$.

In addition, from now on we assume that $S$ is a unital semigroup 
with involution, that is, there exists a
mapping $\fI:S\rightarrow S$ such that $\fI(\fI(a))=a$ and 
$\fI(ab)=\fI(b)\fI(a)$ for all $a,b\in S$. The connection between the
involution $\fI$ and the function $\alpha $ is given by the assumption
\begin{equation}\label{conditie}
\alpha (a\fI (a),x)=1,\quad a\in S,\ x\in X.
\end{equation}

Finally, with the notation and the assumption as before, 
we say that the hermitian kernel $K$ on $X$ is $\phi$-{\em invariant} if
\begin{equation}\label{simetria}
K(x,\phi (a,y))=\overline{\alpha (a,\phi (\fI (a),x))}
\alpha (a,y)K(\phi (\fI (a),x),y)
\end{equation}
for all $x,y\in X$ and $a\in S$. In order to keep the terminology
simple, the function $\alpha$ and the involution $\fI$ will be made
each time precise, if not clear from the context.

The following is the main result of this section.
\begin{theorem}\label{baza}
Let $\phi $ be an action of the unital semigroup 
$S$ with involution $\fI$ satisfying 
\eqref{conditie} on the set $X$ and let $K$ be an $\lh$-valued 
$\phi$-invariant hermitian kernel on $X$. 

The following assertions are equivalent:
\smallskip

{\rm (1)} There exists $L\in {\cB}_{\phi }^+(X,\cH)$ such that $-L\leq K\leq
L$.\smallskip 

{\rm (2)} $K$ has a Kolmogorov decomposition $(V;\cK)$
with the property that there exists a projective
representation $U$ of $S$ on $\cK$ 
(that is, $U(a)U(b)=\sigma (a,b)U(ab)$
for all $a,b\in S$) such that 
\begin{equation}\label{rel}
V(\phi (a,x))=\alpha (a,x)U(a)V(x)
\end{equation}
for all $x\in X$, $a\in S$.
In addition, $\overline{\sigma (\fI (a),a)}U(\fI (a))=
U(a)^{\sharp }$ for all $a\in S$. \smallskip 

{\rm (3)} $K=K_1-K_2$ for two positive definite 
kernels such that $K_1+K_2\in {\cB}_{\phi }^+(X,\cH)$.\smallskip 

{\rm (4)} $K=K_+-K_-$ for two disjoint positive definite 
kernels such that $K_++K_-\in {\cB}_{\phi }^+(X,\cH)$.
\end{theorem}
\begin{proof} (1)$\Rightarrow $(2).
Let ${\cH}_L$ be the Hilbert space
obtained by the completion of the quotient space
${\cF}_0(X,\cH)/{\cN}_L$ with respect to $[\cdot ,\cdot ]_L$, 
where ${\cN}_L=\{f\in {\cF}_0(X,\cH)\mid [f,f]_L=0\}$
is the isotropic subspace of the inner product space
$({\cF}_0(X,\cH),[\cdot ,\cdot ]_L)$. 
Let $A_L$ be the Gram operator
of $K$ with respect to $L$ and let $(V;\cK _{A_L})$ 
be the Kolmogorov decomposition of 
the kernel $K$ described in the previous section.
Since $L$ is $\phi $-bounded, it follows that 
each $\psi _a$ extends to a bounded operator $F(a)$ on ${\cH}_L$.
We notice that 
$$\begin{array}{lll}
[\psi _a(\xi _x),\eta _y]_K&=&
[(\alpha (a,x)^{-1}\xi )_{\phi (a,x)},\eta _y]_K \\[1em]
&=&\alpha (a,x)^{-1}[K(y,\phi (a,x))\xi ,\eta ]_{\cH} \\[1em]
&=&\overline{\alpha (a,\phi (\fI (a),y))} 
[K(\phi (\fI (a),y),x)\xi ,\eta ]_{\cH} \\[1em]
&=&\overline{\alpha (a,\phi (\fI (a),y))}
\overline{\alpha (\fI (a),y)}[\xi _x,\psi _{\fI (a)}(\eta _y)]_K.
\end{array}$$
From the definition of $\sigma $ we have that for $y\in X$,
$$\sigma (a,\fI (a))=\alpha (a,\phi (\fI (a),y))^{-1}
\alpha (\fI (a),y)^{-1}\alpha (a\fI (a),y).$$
By our assumption \eqref{conditie}, $\alpha (a\fI (a),y)=1$, 
so that
$$\sigma (a,\fI (a))=\alpha (a,\phi (\fI (a),y))^{-1}
\alpha (\fI (a),y)^{-1}.$$
Since $|\sigma (a,\fI (a))|=1$, we deduce that
$$[\psi _a(\xi _x),\eta _y]_K=\sigma (a,\fI (a))
[\xi _x,\psi _{\fI (a)}(\eta _y)]_K.$$
This relation can be extended by linearity to 
$$[\psi _a(f),g]_K=\sigma (a,\fI (a))
[f,\psi _{\fI (a)}(g)]_K$$
for all $f,g \in {\cF}_0(X,\cH)$. We deduce that 
$$[A_L\psi _a(f),g]_L=\sigma (a,\fI (a))
[A_Lf,\psi _{\fI (a)}(g)]_L,$$
which implies that 
\begin{equation}\label{comut}
A_LF(a)=\sigma (a,\fI (a))F(\fI (a))^*A_L.
\end{equation}
Theorem~\ref{krein}
implies that there exists a unique operator
$U(a)\in \cL(\cK _{A_L})$
such that 
$$U(a)\Pi _{A_L}=\Pi _{A_L}F(a).$$  
Moreover, for $h \in {\cH}_L$,
$$U(a)U(b)\Pi _{A_L}h=U(a)\Pi _{A_L}F(b)h=\Pi _{A_L}F(a)F(b)h.$$
We also notice that 
$$\begin{array}{lll}
\psi _a\psi _b(\xi _x)&=&\psi _a(\alpha (b,x)^{-1}\xi _{\phi (b,x)}) \\[1em]
&=&\alpha (b,x)^{-1}\alpha (a,\phi (b,x))^{-1}\xi _{\phi (a,\phi
(b,x))} \\[1em]
&=&\sigma (a,b)\alpha (ab,x)^{-1}\xi _{\phi (ab,x)} \\[1em]
&=&\sigma (a,b)\psi _{ab}(\xi _x).
\end{array}$$
We deduce that $F(a)F(b)=\sigma (a,b)F(ab)$ and this relation 
implies that 
$$U(a)U(b)\Pi _{A_L}h=\sigma (a,b)U(ab)\Pi _{A_L}h.$$
Since the set $\{\Pi _{A_L}h\mid h\in \cH _L\}$ is dense 
in $\cK _{A_L}$, we deduce that $U$ is a projective representation
of $S$ on $\cK _{A_L}$.

For $\xi \in \cH$ we have
$$V(\phi (a,x))\xi =\Pi _{A_L}[\xi _{\phi (a,x)}]$$
and 
$$U(a)V(x)\xi =U(a)\Pi _{A_L}[\xi _x]=\Pi _{A_L}F(a)[\xi _x].$$ 
Since $\psi _a(\xi _x)=\alpha (a,x)^{-1}\xi _{\phi (a,x)}$, we 
deduce that
$$F(a)[\xi _x]=\alpha (a,x)^{-1}[\xi _{\phi (a,x)}],$$
so that \eqref{rel} holds.

Finally, the relation \eqref{comut} 
implies that
$$\begin{array}{lll}
[U(a)\Pi _{A_L}f,\Pi _{A_L}g]_{\cK _{A_L}}&=&
[\Pi _{A_L}F(a)f,\Pi _{A_L}g]_{\cK _{A_L}}\\[1em]
 &=&[A_LF(a)f,g]_L \\[1em]
 &=& \sigma (a,\fI (a))[F(\fI (a))^*A_Lf,g]_L \\[1em]
&=&\sigma (a,\fI (a))[A_Lf,F(\fI (a))g]_L \\[1em]
&=&\sigma (a,\fI (a))[\Pi _{A_L}f,\Pi _{A_L}F(\fI (a))g]_{\cK _{A_L}}\\[1em]
&=&\sigma (a,\fI (a))[\Pi _{A_L}f,U(\fI (a))\Pi _{A_L}g]_{\cK _{A_L}}
\end{array}$$
for all $f,g \in \cH _L$, which implies that
$\overline{\sigma (a, \fI (a))}U(\fI (a))=
U(a)^{\sharp }$. We now notice that the
relation \eqref{simetria}
implies that $\sigma (a, \fI (a))=\sigma (\fI (a),a)$, which concludes
the proof of the relation 
$\overline{\sigma (\fI (a),a)}U(\fI (a))=
U(a)^{\sharp }$ for all $a\in S$.

(2)$\Rightarrow $(4).
Let $J$ be a fundamental symmetry on $\cK$. Then 
$J$ is a selfadjoint operator with respect
to the positive definite inner product
$\langle h,g\rangle_J=[Jh,g]_{\cK}$.
Let
$J=J_+-J_-$ be the Jordan decomposition of $J$ and define 
the hermitian kernels 
$$K_{\pm }(x,y)=\pm V(x)^{\sharp }J_{\pm }V(y),\quad 
L(x,y)=V(x)^{\sharp }JV(y), \quad x,y\in X.$$
From $J_++J_-=I$ and $\pm J_{\pm }=J_{\pm }JJ_{\pm }$ we get 
$K(x,y)=K_+(x,y)-K_-(x,y)$ and $L(x,y)=K_+(x,y)+K_-(x,y)$.
To prove that $K_+$ and $K_-$ are positive definite
kernels let $h \in {\cF}_0(X,\cH)$. Then
$$\begin{array}{lll}
\sum\limits_{x,y\in X}
[K_{\pm }(x,y)h(y),h(x)]_{\cH}&=&\sum\limits_{x,y\in X}
[\pm V(x)^{\sharp }J_{\pm }V(y)h(y),h(x)]_{\cH} \\[1.2em]
 &=&\sum\limits_{x,y\in X}
[\pm J_{\pm }V(y)h(y),V(x)h(x)]_{\cK} \\[1.2em]
 &=&\sum\limits_{x,y\in X}
[J_{\pm }JJ_{\pm }V(y)h(y),V(x)h(x)]_{\cK} \\[1.2em]
 &=&\sum\limits_{x,y\in X}
\langle J_{\pm }V(y)h(y),J_{\pm }V(x)h(x)\rangle _J \\[1.2em]
 &=&\|\sum\limits_{x\in X}
J_{\pm }V(x)h(x)\|_J^2\geq 0.
\end{array}$$
It remains to show that $L$ is $\phi $-bounded.
If $h \in{\cF}_0(X,\cH)$, then
$h=\sum _{k=1}^n\xi ^k_{x_k}$ for some $n\in \NN$, vectors
$\xi^1,\ldots,\xi^n\in\cH$ and distinct elements 
$x_1,$ $x_2$, $\ldots $, $x_n$ in $X$. Then

$$\begin{array}{lll}
[\psi _a(h),\psi _a(h)]_L&=&\sum\limits_{j,k=1}^n
[\psi _a(\xi ^j_{x_j}),\psi _a(\xi ^k_{x_k}]_L \\[1.1em]
 &=&\sum\limits_{j,k=1}^n
\alpha (a,x_j)^{-1}
\overline{\alpha (a,x_k)^{-1}}
[\xi ^j_{\phi (a,x_j)},\xi ^k_{\phi (a,x_k)}]_L \\[1.1em]
 &=&\sum\limits_{j,k=1}^n
\alpha (a,x_j)^{-1}
\overline{\alpha (a,x_k)^{-1}}
[L(\phi (a,x_k),\phi (a,x_j))\xi ^j,\xi ^k]_{\cH}\\[1.1em]
 &=&\sum\limits_{j,k=1}^n
\alpha (a,x_j)^{-1}
\overline{\alpha (a,x_k)^{-1}}
\langle V(\phi (a,x_j))\xi ^j,V(\phi (a,x_k))\xi ^k\rangle _J\\[1.1em]
 &=&\sum\limits_{j,k=1}^n
\langle U(a)V(x_j)\xi ^j,U(a)V(x_k)\xi ^k\rangle _J\\[1.1em]
 &=&\|U(a)\sum\limits_{k=1}^nV(x_k)\xi ^k\|_J^2 \\[1.1em]
 &\leq &\|U(a)\|_J^2\|\sum\limits_{k=1}^nV(x_k)\xi ^k\|_J^2\\[1.1em]
&=&\|U(a)\|_J^2\sum\limits_{j,k=1}^n\langle V(x_j)\xi ^j,V(x_k)\xi
^k\rangle _J\\[1.1em]
 &=&\|U(a)\|_J^2\sum\limits_{j,k=1}^n[\xi ^j_{x_j},\xi ^k_{x_k}]_L \\[1.1em]
 &=&\|U(a)\|^2_J[h,h]_L,
\end{array}$$
\noindent
so that $L$ is $\phi $-bounded.

We also deduce that $(V,(\cK ,\langle \cdot ,\cdot \rangle _J))$
is the Kolmogorov decomposition of the 
positive definite kernel $L$
and $(J_{\pm }V,(J_{\pm }\cK ,\langle \cdot ,\cdot \rangle _J))$
is the Kolmogorov decomposition of $K_{\pm }$.
Since $J_+J_-=0$ we deduce that $J_+\cK\cap J_-\cK=\{0\}$ and, 
by Proposition~16, in \cite{Sch} we deduce that $K_+$ and $K_-$ are
disjoint kernels. 

Since (4)$\Rightarrow $(3)
and (3)$\Rightarrow $(1) are obvious implications, the proof is complete.
\end{proof}

A Kolmogorov decomposition $(V,\cK )$ of the hermitian 
kernel $K$ 
for which there exists a projective representation
$U$ such that \eqref{rel} holds is called a
{\it projectively invariant Kolmogorov decomposition}.
Also, a projective representation $U$ satisfying the additional
property $U(a)^{\sharp }=\overline{\sigma (\fI (a),a)}U(\fI (a))$ 
for all $a\in S$, is called {\it symmetric projective representation}. 

A natural question that can be raised in connection with the previous
result is whether ${\cB}_{\phi }^+(X,\cH)$ is a sufficiently rich
class of kernels. 

\begin{proposition}\label{invizometric}
Assume that $S$ is a group and $\fI(a)=a^{-1}$, $a\in S$. If $K$ is a
$\phi$-invariant hermitian kernel on $X$ then, for any $a\in S$ the
operator $\psi_a$ is isometric
with respect to the inner product $[\cdot,\cdot]_K$. In particular,
any $\phi$-invariant positive definite kernel on $X$ belongs to 
${\cB}_{\phi }^+(X,\cH)$.
\end{proposition}

\begin{proof} Indeed, in this case \eqref{conditie} becomes
$\alpha(e,x)=1$ for all $x\in X$, where $e$ is the unit of the group $S$.
Also, if $K$ is a hermitian kernel then it is $\phi$-invariant if and
only if
$$K(\phi(a,x),\phi (a,y))=\overline{\alpha (a,x)}
\alpha (a,y)K(x,y),\quad x,y\in X,\ a\in S.$$
Let $\xi ,\eta \in \cH$ be arbitrary. Then 
$$\begin{array}{lll}
[\psi _a(\xi _x),\psi _a(\eta _y)]_K&=&
\alpha (a,x)^{-1}\overline{\alpha (a,y)}^{-1}
[\xi _{\phi (a,x)},\eta _{\phi (a,y)}]_K \\[1em]
&=&
\alpha (a,x)^{-1}\overline{\alpha (a,y)}^{-1}
[K(\phi (a,y),\phi (a,x))\xi ,\eta ]_{\cH} \\[1em]
&=&[K(y,x)\xi ,\eta ]_{\cH}=[\xi _x,\eta _y]_K,
\end{array}$$
and hence $\psi_a$ is $[\cdot,\cdot]_K$ isometric. \end{proof}

\begin{remark} {\em
$(1)$ Theorem \ref{baza} is known when 
$\cH$ is a Hilbert space and the kernel $K$ is 
positive definite and satisfies
\begin{equation}\label{uu1}
K(\phi (a,x),\phi (a,y))=\overline{\alpha (a,x)}\alpha
(a,y)K(x,y),\quad a\in S,\ x,y\in X.
\end{equation}
(see, for instance, \cite{PS}).
In that case the proof is easily obtained by defining directly
\begin{equation}\label{uu2}
U(a)V(x)\xi =\alpha (a,x)^{-1}V(\phi (a,x))\xi 
\end{equation}
for $\xi \in \cH$ and verify that $U(a)$ satisfies all the 
required properties (we note that no involution
is considered in this case). We have to emphasize that this direct 
approach does not work in the hermitian case since the 
formula \eqref{uu2}
does not necessarily give a bounded operator.
In order to overcome this difficulty we have to replace
the symmetry condition in \eqref{uu1} by the symmetry condition in 
\eqref{simetria}
and then use Theorem~\ref{krein}. This was the main point 
in the proof of Theorem~\ref{baza}.

$(2)$ The positive definite version of 
Theorem~\ref{baza} has many applications, 
some of them mentioned for instance in \cite{EL}, \cite{EK}, 
and \cite{PS}. Such a typical application gives a Naimark
dilation for Toeplitz kernels. Thus, if $X=S$,
$\phi (a,x)=ax$, and $\alpha (a,x)=1$
for all $a,x\in S$, then \eqref{uu1} becomes the well-known
Toeplitz condition
$$K(ab,ac)=K(b,c)$$
for all $a,b,c\in S$. If $K$ is a positive definite kernel
on $S$ satisfying the Toeplitz condition and $K(e,e)=I$, 
where $e$ is the unit of $S$, then $\{U(a)\}_{a\in S}$
defined by \eqref{uu2} is a semigroup of isometries on a
Hilbert space $\cK$ containing $\cH$ such that
$$K(a,b)=P_{\cH}U(a)^*U(b)|\cH ,$$
for all $a,b\in S$, where $P_{\cH }$ denotes the 
orthogonal projection of $\cK$ onto $\cH$.

$(3)$ The next example explores the fact that for positive definite kernels
the representation $\{U(a)\}_{a\in S}$ given by
\eqref{uu2} is unique up to unitary equivalence.
Thus, consider the action of a group $G$ on the 
Hilbert space $\cH$ such that 
$\langle \phi (g,\xi ),\phi (g,\eta )\rangle 
=\langle \xi ,\eta \rangle $
for all $g\in G$ and $\xi , \eta \in \cH$. We consider the 
kernel 
$K(\xi , \eta )=\langle \eta ,\xi \rangle $
on $\cH $ and notice that $K$ is positive definite.
Its Kolmogorov decomposition is given by
$V(\xi ):{\CC }\rightarrow \cH$, 
$V(\xi )\lambda =\lambda \xi ,$
$\lambda \in \CC ,\xi \in \cH $.
If we use the positive definite version of Theorem~\ref{baza}, 
we deduce that there exists a 
Kolmogorov decomposition $V'$ of $K$ and a representation
$U'$ of $G$ such that 
$V'(\phi (g,\xi ))=U'(g)V'(\xi )$
for all $g\in G$ and $\xi \in \cH $. 
From the uniqueness of $V$ up to unitary equivalence, it
follows that there exists a unitary operator $\Phi $
such that 
$V(\phi (g,\xi ))=\Phi U'(g) \Phi ^*V(\xi ),$
or
$\phi (g,\xi )=U(g)\xi ,$
with $U(g)=\Phi U'(g) \Phi ^*$. 
Therefore we obtained the well-known result that 
$\phi $ acts by linear unitary operators.}
\end{remark}

The last example was intended to emphasize the importance of  
the uniqueness up to unitary equivalence of the projectively invariant 
Kolmogorov decompositions. This issue turns out to be 
rather delicate in the hermitian case. Theorem~\ref{unic}
settles this question only in the case of the trivial 
semigroup $S$ with one element. 
It is easily seen that the spectral condition in Theorem~\ref{unic}
is also sufficient for the uniqueness of a projectively invariant 
Kolmogorov decomposition. However, Theorem~\ref{contraexemplu}
shows that the proof in \cite{CG1} of Theorem~\ref{unic}
cannot be easily adapted to the case of an arbitrary semigroup $S$.

We conclude this section with another case when uniqueness holds.
Given a hermitian kernel $K$, the rank $\rank(K)$ is, by definition, the
supremum of $\rank(K_\Delta )$ taken over all finite subsets
$\Delta \subset X$, where $K_\Delta $ is the
restricted kernel $(K(x,y))_{x,y\in\Delta }$. By definition $\rank(K)$ is
either a positive integer or the symbol $\infty$.
A hermitian kernel $K$ has $\kappa$ {\it negative squares} 
if the inner product space $(\cF _0(X,\cH), [\cdot,\cdot]_K)$ has 
negative signature
$\kappa$, that is, $\kappa $ is the maximal dimension of all its negative
subspaces. It is easy to see that this is equivalent to $K=K_+-K_-$,
where $K_\pm $ are disjoint positive definite kernels such that
$\rank(K_-)=\kappa $, see e.g.\ \cite{Sch}. This allows us to
define $\kappa ^-(K)=\kappa $, the {\em number of negative squares} of the
kernel $K$. In particular, hermitian kernels with a finite
number of negative squares always have Kolmogorov decomposition and
for any Kolmogorov decomposition $(V;\cK)$ of $K$ we have
$\kappa^-(\cK)=\kappa^-(K)<\infty$, hence $\cK$ is a 
Pontryagin space with negative signature
$\kappa$. 

In Pontryagin spaces the strong topology is intrinsicly
caracterized in terms of the indefinite inner product, e.g.\ see
\cite{IKL}. Therefore, by using Proposition \ref{invizometric} and 
Shmul'yan's Theorem
(e.g.\ see Theorem~2.10 in \cite{DR}) we get:

\begin{theorem}\label{pontr}
Let $\phi $ be an action of the group 
$S$ on the set $X$ and let $K$ be an $\lh$-valued 
$\phi$-invariant hermitian kernel on $X$ with a finite number 
of negative squares. Then $K$ has a projectively invariant 
Kolmogorov decomposition on a Pontryagin space, that is 
unique up to unitary equivalence.
\end{theorem}

\section{Similarity}

The symmetric projective representation $U$ of $S$ obtained
in Theorem~\ref{baza} acts on a Kre\u\i n space. It would be of
special interest to decide whether $U$ is at least similar to a symmetric
projective representation on a Hilbert space, a property
related to the well-known similarity problem for group 
representations, see \cite{Pis}
for a recent survey.

The above mentioned problem is also closely related to the
characterization of 
those $\phi$-invariant hermitian kernels $K$ with the property that 
the representation $K=K_+-K_-$
holds for two positive definite $\phi$-invariant kernels.

In this section we give an answer to these two questions in terms of
fundamental reducibility.
We say that the projective representation 
$U$ of $S$ on the Krein space $\cK$ is 
{\it fundamentally
reducible} if there exists a fundamental 
symmetry $J$ on $\cK$ such that $U(a)J=
JU(a)$ for all $a\in S$. This condition is readily equivalent to the
condition $U(a)^\sharp=U(a)^*$ for all $a\in S$, and further,
equivalent to the diagonal representation of $U(a)$ with respect to
a fundamental decomposition of the Kre\u\i n space $\cK$.

\begin{proposition}\label{invsimilaritate} Let $S$ be a semigroup with
involution $\fI$ and $\sigma$ satisfies the $2$-cocycle 
property \eqref{cociclu} on $S$. Let $U$ be a symmetric projective 
representation of $S$ on the Kre\u\i n space $\cK$.
Then the following assertions are equivalent:

{\rm (1)} $U$ is similar to a symmetric projective representation
$T$ on a Hilbert space.

{\rm (2)} $U$ is fundamentally reducible.
\end{proposition}

\begin{proof} (1)$\Rightarrow $(2). Let $\Phi\in\cL(\cK,\cG)$
be the similarity such that $T(a)\Phi=\Phi U(a)$ for $a\in S$. We
first notice that $\Phi $ is also an {\it involutory similarity}
(with the
terminology from \cite{KS}), that is
\begin{equation}\label{invsim}T(a)^*=\Phi U(a)^\sharp \Phi^{-1},
\quad a\in S.\end{equation}
Then, we consider on $\cK$ the positive inner product
$\langle\xi,\eta\rangle_\Phi=\langle \Phi\xi,\Phi\eta\rangle$,
$\xi,\eta\in\cK$. Since $\Phi$ is boundedly invertible, there exists a
selfadjoint and boundedly invertible operator $G\in\cL(\cK)$ such that
$[\xi,\eta]=\langle G\xi,\eta\rangle_\Phi$,
$\xi,\eta\in\cK$. Therefore, for arbitrary $a\in S$ and
$\xi,\eta\in\cK$ we have
$$\begin{array}{rcl} \langle U(a)\xi,\eta\rangle_{\Phi} & = & \langle
\Phi U(a)\xi,\Phi\eta\rangle=\langle T(a)\Phi\xi,\Phi\eta\rangle \\[1em]
& = & \langle \Phi\xi,T(a)^*\Phi\eta\rangle=\langle \Phi\xi,\Phi
U(a)^\sharp \eta\rangle \\[1em]
& = & \langle\xi,U(a)^\sharp\eta\rangle_\Phi
=[G^{-1}\xi,U(a)^\sharp\eta] \\[1em]
& = & [U(a)G^{-1}\xi,\eta]=
\langle GU(a)G^{-1}\xi,\eta\rangle_\Phi.\end{array}$$ 
Thus, $GU(a)=U(a)G$ and letting $J=\sgn(G)$ it follows that $J$ is a
fundamental symmetry on the Kre\u\i n space $\cK$ such that
$JU(a)=U(a)J$.\smallskip

(2)$\Rightarrow $(1). If $J$ is a
fundamental symmetry on the Kre\u\i n space $\cK$ such that
$JU(a)=U(a)J$, for all $a\in S$, 
then $U$ is a symmetric projective representation with
respect to the Hilbert space $(\cK,\langle\cdot,\cdot\rangle_J)$.\end{proof}

With the notation as in Proposition \ref{invsimilaritate}, if
$\sigma$ has the $2$-cocycle property \eqref{cociclu}
and $|\sigma(a,b)|=1$ for all
$a,b\in S$, then it follows that 
\begin{equation}\label{izometric}U(a)^\sharp U(a)=U(\fI(a)a), \quad
a\in S.\end{equation} 
Thus, in certain applications where $U$ consists of (Kre\u\i n
space) isometric operators,
it is interesting to know whether $U$ is similar to a symmetric projective
representation of isometric
operators on a Hilbert space. Clearly, a necessary condition is that for some
(equivalently for all) unitary norm $\|\cdot\|$ on $\cK$ there exists
$C>0$ such that
\begin{equation}\label{marginire} \frac{1}{C}\|\xi\|\leq
\|U(a)\xi\|\leq C\|\xi\|,\quad a\in S,\ \xi\in\cK.\end{equation}
As expected, the converse implication is related to the assumption of
amenability of the semigroup $S$. More precisely, following closely the
idea in the proof of Th\`eor\'eme 6 in \cite{Dix}, we get:

\begin{theorem}\label{dixmier} Let $S$ be an amenable semigroup,
$\sigma$ has the $2$-cocycle property \eqref{cociclu}, 
$|\sigma(a,b)|=1$ for all
$a,b\in S$, and let $U$ be a projective representation (without any
assumption of symmetry) of
$S$ on a Hilbert space $\cK$, such that \eqref{marginire} holds for
some constant $C>0$. Then $U$ is similar to a projective
representation $T$ of $S$ on a Hilbert space $\cG$
such that $T(a)$ are isometric for all $a\in S$.\end{theorem}

We come now to the problem of characterizing those hermitian invariant
kernels that can be represented as a difference of two positive
invariant kernels.

\begin{theorem}\label{balabusta}
Let $\phi $ be an action of the unital semigroup
$S$ with involution $\fI$ satisfying
\eqref{conditie} on the set $X$ and let $K$ be an $\lh$-valued 
$\phi $-invariant hermitian kernel on $X$.
The following assertions are equivalent:

{\rm (1)} There exists $L\in {\cB}_{\phi }^+(X,\cH)$ such that $-L\leq K\leq
L$ and $L$ is $\phi$-invariant.

{\rm (2)} $K$ has a projectively invariant
Kolmogorov decomposition $(V;\cK)$ such that 
the associated projective
representation is fundamentally reducible.

{\rm (3)} $K=K_+-K_-$ for two disjoint positive definite
kernels such that $K_++K_-\in {\cB}_{\phi }^+(X,\cH)$ and both $K_\pm$
are $\phi$-invariant.
\end{theorem}
\begin{proof} (1)$\Rightarrow $(2).
We use the same notation as in the proof of Theorem~\ref{baza}.
Thus, ${\cH}_L$ is the Hilbert space
obtained by the completion of the quotient space
${\cF}_0(X,\cH)/{\cN}_L$ with respect to $[\cdot ,\cdot ]_L$,
where ${\cN}_L$
is the isotropic subspace of the inner product space
$({\cF}_0(X,\cH),[\cdot ,\cdot ]_L)$.
Let $A_L$ be the Gram operator
of $K$ with respect to $L$ and let $(V;\cK _{A_L})$
be the projectively invariant Kolmogorov decomposition of
the kernel $K$ described in the proof of (1)$\Rightarrow $(2)
in Theorem~\ref{baza}.
Since $L$ is $\phi $-bounded, it follows that
each $\psi _a$ extends to a bounded operator $F(a)$ on ${\cH}_L$.
Since $L$ is $\phi$-invariant, we deduce that
$$[\psi _a(f),g]_L=\sigma (a,\fI (a))
[f,\psi _{\fI (a)}(g)]_L$$
for all $f,g \in {\cF}_0(X,\cH)$, which implies that
$$F(a)=\sigma (a,\fI (a))F(\fI (a))^*.$$
This relation and \eqref{comut} imply that 
$$A_LF(a)=F(a)A_L$$
for all $a\in S$.
Let $A_L=S_{A_L}|A_L|$ be the polar
decomposition of $A_L$ and let $J_{A_L}$ be
the symmetry introduced in Example~\ref{hindus}.
Using \eqref{JPI}, we deduce that
$$\begin{array}{lll}
U(a)J_{A_L}\Pi _{A_L}&=&U(a)\Pi _{A_L}S_{A_L} \\[1em]
&=&\Pi _{A_L}F(a)S_{A_L} \\[1em]
&=&\Pi _{A_L}S_{A_L}F(a) \\[1em]
&=&J_{A_L}\Pi _{A_L}F(a) \\[1em]
&=&J_{A_L}U(a)\Pi _{A_L},
\end{array}$$
therefore the representation $U$ is fundamentally reducible.

(2)$\Rightarrow $(3).
We consider the elements involved in the
proof of (2)$\Rightarrow $(4) in Theorem~\ref{baza}
for a fundamental symmetry $J$ on $\cK$ for which
$U(a)J=JU(a)$, $a\in S$.
Therefore $U(a)J_{\pm }=J_{\pm }U(a)$ for all $a\in S$, and then
$$\begin{array}{lll}
K_{\pm }(x,\phi (a,y))&=&\pm V(x)^{\sharp }J_{\pm }V(\phi (a,y)) \\[1em]
&=& \pm \alpha (a,y)V(x)^{\sharp }J_{\pm }U(a)V(y) \\[1em]
&=& \pm \alpha (a,y)V(x)^{\sharp }U(a)J_{\pm }V(y) \\[1em]
&=& \pm \alpha (a,y)\sigma (\fI (a),a)V(x)^{\sharp }U(\fI (a))^{\sharp }
J_{\pm }V(y) \\[1em]
&=& \pm \overline{\alpha (a,\phi (\fI (a),x))}
\alpha (a,y)V(\phi (\fI (a),x))^{\sharp }J_{\pm }V(y) \\[1em]
&=& \overline{\alpha (a,\phi (\fI (a),x))}
\alpha (a,y)K_{\pm }(\phi (\fI (a),x),y)).
\end{array}$$

(3)$\Rightarrow $(1).
Just set $L(x,y)=K_+(x,y)+K_-(x,y)$.
\end{proof}

In case $S$ is a group with the involution $\fI(a)=a^{-1}$, then
some of the assumptions in the previous results
simplify to a certain extent. In this case, as a consequence of
\eqref{izometric}, the symmetric projective representation $U$ associated to a
$\phi$-invariant Kolmogorov decomposition consists of unitary operators.

\begin{theorem}\label{similaritate} Let $S$ be a group
and $\sigma$ a $2$-cocycle on $S$ with $|\sigma(a,b)|=1$ for all
$a,b\in S$. Let $U$ be a unitary
projective representation of $S$ on the Kre\u\i n space
$\cK$. Then the following assertions are equivalent:

{\rm (1)} $U$ is similar to a unitary projective representation
$T$ on a Hilbert space, that is, $T\colon S\rightarrow\cL(\cG)$, $\cG$
a Hilbert space, $T(a)T(b)=\sigma(a,b)T(ab)$ and
$T(a)^* =\sigma(a^{-1},a)T(a^{-1})$ for all $a\in S$.

{\rm (2)} $U$ is fundamentally reducible.

Moreover, if $U$ satisfies one (hence both) of the assumptions {\rm (1)} and 
{\rm (2)} then $U$ is uniformly bounded, that is,
\begin{equation}\label{marginire2} \sup_{a\in
S}\|U(a)\|<\infty.\end{equation}
If, in addition, $S$ is amenable, then \eqref{marginire2} is equivalent to
(any of) the conditions {\rm (1)} and {\rm (2)}.\end{theorem}

\begin{proof} This follows from Proposition \ref{invsimilaritate} and
Theorem \ref{dixmier}.\end{proof}

\begin{theorem}\label{corolar}
Let $\phi $ be an action of the group
$S$ on the set $X$ and let $K$ be an $\lh$-valued $\phi$-invariant
hermitian kernel on $X$.

The following assertions are equivalent:
\smallskip

{\rm (1)} There exists a $\phi$-invariant positive definite  
$L$ on $X$ such that $-L\leq K\leq L$.

{\rm (2)} $K$ has a projectively invariant
Kolmogorov decomposition $(V;\cK)$ such that
the associated symmetric projective representation is similar to a symmetric
projective representation on a Hilbert space.

{\rm (3)} $K=K_+-K_-$ for two disjoint positive definite $\phi$-invariant
kernels.\end{theorem}
\begin{proof}
This follows from Proposition \ref{invizometric} and Theorem~\ref{balabusta}.
\end{proof}

\section{An example: Weyl exponentials}

In this section we discuss an example leading to 
a projectively invariant hermitian kernel.
Thus, the Fock representation of the canonical 
commutation relations is obtained from an action 
of the rigid motions of a Hilbert space on the 
exponential vectors of a Fock space. It is natural to consider 
a similar construction for other groups, like Poincar\'e group,
involving indefinite inner products. Various models involving
Fock spaces associated to indefinite inner products
were studied in \cite{MPS}, \cite{Str}.
Here we emphasize that the Kolmogorov decomposition gives a simple 
construction of the Weyl exponentials (the related topic of the 
representations of the Heisenberg algebra in Kre\u{\i}n spaces 
is taken up in 
\cite{MMSV}).

Let $(\cH ,[\cdot ,\cdot ])$ be a Kre\u{\i}n space and 
consider $\cP$ the group of its rigid motions. This is
the semidirect product of the additive group $\cH$
and the group of the bounded unitary operators on $\cH$. 
The group law is given by
$$(\xi ,U)(\xi ',U')=(\xi+U\xi ',UU')$$
and an action of $\cP$ on $\cH$ can be defined by the formula
$$\phi ((\xi ,U),\xi ')=\xi+U\xi '.$$
In particular, the normal subgroup $\cH$ of $\cP$ acts on $\cH$ 
by translations. For simplicity, we restrict here 
to this action by translations. The hermitian kernel associated to 
this construction is 
defined by the formula:
\begin{equation}\label{Fock1}
K(\xi ,\eta )=\exp (\frac{i\Im [\eta ,\xi ]}{2})
\exp (-\frac{[\xi -\eta ,\xi -\eta ]}{4}),
\end{equation}
for $\xi ,\eta \in \cH $.
The additive group $\cH$ acts on itself by the 
translations  $\phi (\xi ,\eta )=\xi +\eta $ and we notice that
\begin{equation}\label{Fock2}
K(\phi (\xi ,\eta ),\phi (\xi ,\eta '))=\overline{\alpha (\xi ,\eta )}
\alpha (\xi ,\eta ')K(\eta ,\eta ')
\end{equation}
for all $\xi ,\eta ,\eta '\in \cH$, where
$$\alpha (\xi ,\eta )=\exp (-\frac{i\Im [\xi ,\eta ]}{2})$$
and then 
$$\sigma (\xi ,\eta )=
\alpha (\xi ,\eta +\eta ')^{-1}\alpha (\eta ,\eta ')^{-1}
\alpha (\xi +\eta ,\eta ')=\exp (\frac{i\Im [\xi ,\eta ]}{2}).$$

\begin{theorem}\label{Weylexp}
The kernel $K$ defined by \eqref{Fock1} has a Kolmogorov
decomposition $(V,\cK )$ with the property that
the operators defined by the formula
\begin{equation}\label{Fock3}
\alpha (\xi ,\eta )W(\xi )V(\eta )=V(\xi +\eta )
\end{equation}
are defined on the common dense domain
$\bigvee _{\xi \in \cH}V(\xi )\CC $ in $\cK$ 
and satisfy the canonical commutation relations
\begin{equation}\label{Fock4}
W(\xi )W(\eta )=\sigma (\xi ,\eta )W(\xi +\eta ).
\end{equation} 
\end{theorem}

\begin{proof}
We can obtain a Kolmogorov decomposition of the kernel
$K$ by adapting the Fock space construction from the 
positive definite case.
Let $J$ be a fundamental symmetry on $\cH $ and 
let $(\cH ,\langle \cdot ,\cdot \rangle _J)$
be the associated Hilbert space structure on $\cH $.
We then consider the Fock space
$$F(\cH )=\bigoplus\limits _{n=0}^{\infty }\cH _n,$$
where $\cH _0$ is a one-dimensional Hilbert space
generated by the unit vector $\Omega $
and $\cH _n$ is the $n$-fold tensor product
$\otimes ^n\cH $. Then,
the operator
$$J_F=\bigoplus\limits _{n=0}^{\infty }\otimes ^nJ$$
is a symmetry of $F(\cH )$. Let $P_n=(n!)^{-1}\sum _{\pi \in S_n}\pi $
be the orthogonal projection of $\cH _n$ onto its
symmetric part, where
$$\pi (\xi _1\otimes \ldots \otimes \xi _n)=\xi _{\pi ^{-1}(1)}
\otimes \ldots \otimes \xi _{\pi ^{-1}(n)}$$
and $\pi $ is an element of the permutation group $S_n$
on $n$ symbols. We notice that
$$P_n(\otimes ^nJ)=(\otimes ^nJ)P_n$$
for all $n\geq 0$, therefore
$$(\bigoplus\limits_{n=0}^{\infty}P_n)J_F=
J_F(\bigoplus\limits_{n=0}^{\infty}P_n).$$
It follows that the compression of $J_F$ to 
$F^s(\cH )=(\oplus _{n=0}^{\infty }P_n)F(\cH )$, the boson Fock space, 
gives a symmetry on $F^s(\cH )$. So, $F^s(\cH )$
becomes a Kre\u{\i}n space by setting $[x,y]_{F^s(\cH )}=
\langle J_Fx,y \rangle $, $x,y\in F^s(\cH )$ and 
$\langle \cdot ,\cdot \rangle $ is the positive definite
inner product on $F^s(\cH )$ obtained from 
$\langle \cdot ,\cdot \rangle _J$. We now use the Kolmogorov
decomposition that gives the Bose-Fock space (see
\cite{EL} or \cite{PS} for more details). Thus, for
$\xi \in \cH $ we define
the map
$$\Exp(\xi )=\bigoplus\limits _{n=0}^{\infty }(n!)^{-1/2}\xi _n,$$
where $\xi _0=\Omega $ and $\xi _n=\xi \otimes \ldots \otimes \xi$.
Then, we define
\begin{equation}\label{Fock5}
V(\xi )\lambda =\lambda 
\exp (-\frac{[\xi ,\xi ]}{4})\Exp(\frac{\xi }{\sqrt{2}})
\end{equation}
and we can easily check that $(V,F^s(\cH ))$
is a Kolmogorov decomposition of $K$.
We now define $\cD (W)=\bigvee _{\xi \in \cH}V(\xi ){\CC}$
and each element of $\cD (W)$ admits a representation
of the form $\sum _{k=1}^nz_kV(\xi _k)\lambda _k$,
where $z_1,$ $\ldots $, $z_n$, 
$\lambda _1$, $\ldots $, $\lambda _n$ are complex numbers and 
$\xi _1 $, $\ldots $, $\xi _n$ are elements in $\cH$.
For each $\xi \in \cH$ we can define a map on 
$\cD (W)$ by the formula
$$W(\xi )(\sum _{k=1}^nz_kV(\xi _k)\lambda _k)=
\sum _{k=1}^nz_k\alpha (\xi ,\xi _k)^{-1}V(\xi +\xi _k)\lambda _k.$$
In order to show that this map is well-defined
it is enough to prove that the relation 
$\sum _{k=1}^nz_kV(\xi _k)\lambda _k=0$
implies $\sum _{k=1}^nz_k\alpha (\xi ,\xi _k)^{-1}V(\xi +\xi _k)\lambda _k=0$.
This follows from the fact that $\{\Exp(\xi )\mid \xi \in \cH \}$
is a linearly independent set in $F^s(\cH )$ (Theorem~6.3 in \cite{EL}).
Next we have
$$\begin{array}{lll}
W(\xi )W(\eta )V(\eta ')&=&W(\xi )\alpha (\eta ,\eta ')^{-1}V(\eta
+\eta ') \\[1em]
&=&\alpha (\eta ,\eta ')^{-1}\alpha (\xi ,\eta +\eta ')^{-1}
V(\xi +\eta +\eta ') \\[1em]
&=&\sigma (\xi ,\eta )\alpha (\xi +\eta ,\eta ')^{-1}V(\xi +\eta +\eta
') \\[1em] 
&=&\sigma (\xi ,\eta )W(\xi +\eta )V(\eta '),
\end{array}$$
so that the family $\{W(\xi )\}_{\xi \in \cH }$ satisfies the 
canonical commutation relations. We also notice that the
operators $W(\xi )$ are isometric on the common domain $\cD (W)$
with respect to the indefinite inner product on $F^s(\cH )$.
\end{proof}

\section{Representations of $*$-algebras}

In this section we discuss the GNS representation for unital 
$*$-algebras from the point of view of invariant hermitian 
kernels. Our goal is to make connections with some 
constructions of interest in quantum field theories
such as those summarized in \cite{Str}.

Let $\cA$ be a $*$-algebra with identity $1$ and let
$Z$ be a linear hermitian functional on $\cA$ with 
mass $1$ ($Z(1)=1$). Then $\cA$ is a unital multiplicative
semigroup with involution acting on itself by
\begin{equation}\label{actstar}
\phi (a,x)=xa^*
\end{equation}
for $a,x\in \cA$. We define
\begin{equation}\label{nucstar}
K_Z(x,y)=Z(xy^*)
\end{equation}
for $x,y\in \cA$. Then $K_Z$ is a hermitian
kernel on $\cA$ with scalar values and satisfies the
symmetry relation
\begin{equation}\label{simstar}
K_Z(x,\phi (a,y))=Z(xay^*)=K_Z(\phi (a^*,x),y)
\end{equation}
for $a,x,y\in \cA$. 
In order to describe the GNS construction for $Z$ we will 
use the concept of unbounded representations of $\cA$.
Thus, a mapping $\pi $ of $\cA$ into the set of 
closable operators defined on a common dense domain 
$\cD (\pi )$ of a Banach space
$\cK$ is called a {\it closable representation} if it is linear, $\cD(\pi)$ 
is invariant under all operators $\pi(a)$, $a\in\cA$, and 
$\pi(ab)=\pi(a)\pi(b)$ for all $a,b\in \cA$. If, in addition, 
$\cK$ is a Kre\u\i n space and, for all $a\in\cA$, 
the domain of  $\pi (a)^{\sharp}$
contains $\cD (\pi )$ and 
\begin{equation}\label{adj}
\pi (a)^{\sharp}|\cD (\pi )=\pi (a^*),
\end{equation}
then $\pi$ is called a {\it hermitian closable representation} on the
Kre\u\i n  space $\cK$ (or, a $J$-{\it representation}, as introduced
in \cite{ota}, see also \cite{hof}). 

The {\it GNS data} 
$(\pi ,\cK ,\Omega)$ associated to $Z$
consists of a hermitian closable representation of $\cA$ on the 
Kre\u{\i}n space $\cK$ and a vector $\Omega\in \cD (\pi )$ such that
\begin{equation}\label{gns}
Z(a)=[\pi (a)\Omega ,\Omega]_{\cK}
\end{equation}
for all $a\in \cA$ and $\bigvee _{a\in \cA}\pi (a)\Omega=\cD (\pi )$.
It was known that not every hermitian functional $Z$ admits GNS data.
Characterizations of those $Z$ that do admit GNS data
appeared in papers such as \cite{MPS}, \cite{AGW}, \cite{hof}. 
We show that the GNS data associated to a hermitian form can be
equivalently described in terms of Kolmogorov decompositions of 
the kernel $K_Z$.

\begin{proposition}\label{GNSkol}
Let $\cA$ be a unital $*$-algebra, let $Z$ be a linear
hermitian functional on $\cA$ with $Z(1)=1$, and consider
the kernel $K_Z$ associated to $Z$ by \eqref{nucstar}. For every GNS
data $(\pi,\cK,\Omega)$ of $Z$ define
\begin{equation}\label{kolGNS} V(a)\lambda=\pi(a^*)\lambda\Omega,\quad
a\in\cA,\ \lambda\in\CC.\end{equation}
Then $(V,\cK)$ is a Kolmogorov decomposition of the hermitian kernel
$K_Z$ and \eqref{kolGNS} establishes a bijective correspondence
between the set of all GNS data of $Z$ and the set of all Kolmogorov
decompositions of $K_Z$.

In particular, $Z$ admits GNS data if and only if the hermitian hernel
$K_Z$ has Kolmogorov decompositions.
\end{proposition}

\begin{proof} Let $(\pi,\cK,\Omega)$ be a GNS data of $Z$. Then for
arbitrary $a\in\cA$, $V(a)$ defined as in \eqref{kolGNS} is a linear
bounded mapping from $\CC$ into $\cK$ such that 
$\bigvee_{a\in cA}V(a)\CC$ is dense in $\cK$. Also, 
$$\begin{array}{lll}
V(x)^{\sharp }V(y)&=&[V(y)1,V(x)1]_{\cK}\\[1em]
 &=&[\pi (y^*)\Omega,\pi (x^*)\Omega]_{\cK}\\[1em]
 &=&[\pi (x)\pi (y^*)\Omega,\Omega]_{\cK}\\[1em]
 &=&[\pi (xy^*)\Omega,\Omega]_{\cK}\\[1em]
 &=&Z(xy^*)=K_Z(x,y),
\end{array}$$
and hence $(V,\cK )$ is a Kolmogorov decomposition of $K_Z$.

Let $(V,\cK)$ be a Kolmogorov decomposition
of $K_Z$. Define $\cD (\pi )=\bigvee _{x\in \cA}V(x)\CC$
which is dense in $\cK$. Each element of $\cD (\pi )$ admits a 
representation of the form $\sum _{k=1}^n\alpha _kV(x_k)\lambda _k$, 
where $\alpha _1$, $\ldots $, $\alpha _n$,
$\lambda _1$, $\ldots $, $\lambda _n\in {\CC}$,
$x_1$, $\ldots $, $x_n\in \cA$.
We can define a map on $\cD (\pi )$ by the formula
\begin{equation}\label{pia}\pi (a)(\sum _{k=1}^n\alpha _kV(x_k)\lambda _k)=
\sum _{k=1}^n\alpha _kV(x_k a^*)\lambda _k.\end{equation}
We have to prove that this map is well-defined. To this end, it is enough
to show that $\sum _{k=1}^n\alpha _kV(x_k)\lambda _k=0$
implies $\sum _{k=1}^n\alpha _kV(x_ka^*)\lambda _k=0.$
From $\sum _{k=1}^n\alpha _kV(x_k)\lambda _k=0$
we deduce that for any $y\in \cA$ we have
$$\begin{array}{lll}
0 &=& V(y)^{\sharp }(\sum\limits _{k=1}^n\alpha _kV(x_k)\lambda _k) \\[1em]
 &=&\sum\limits _{k=1}^n\alpha _kV(y)^{\sharp }V(x_k)\lambda _k \\[1em]
 &=&\sum\limits _{k=1}^n\alpha _kK(y,x_k)\lambda _k.
\end{array}$$
On the other hand, for an arbitrary $z\in \cA$, we have
$$\begin{array}{lll}
V(z)^{\sharp }(\sum\limits _{k=1}^n\alpha _kV(x_k a^*)\lambda _k)
&=&\sum\limits _{k=1}^n\alpha _kV(z)^{\sharp }V(x_k a^*)\lambda _k \\[1em]
 &=&\sum\limits _{k=1}^n\alpha _kK(z,x_k a^*)\lambda _k \\[1em]
 &=&\sum\limits _{k=1}^n\alpha _kK(za,x_k)\lambda _k,
\end{array}$$
whence letting $y=za$, the last sum is zero. 
This implies that for arbitrary $\lambda \in \CC$,
$$[\sum _{k=1}^n\alpha _kV(x_ka^*)\lambda _k,V(y)\lambda ]_{\cK}=0.$$
But $\bigvee _{x\in \cA}V(x)\CC $ is dense in $\cK$, therefore
$\sum _{k=1}^n\alpha _kV(x_ka^*)\lambda _k=0$ and hence
$\pi (a)$ is well-defined on $\cD (\pi )$. Also, $\pi (a)$ is 
a linear operator on $\cD (\pi )$. We easily check that
$\pi (a)\pi (b)f=\pi (ab)f$ for $f\in \cD (\pi )$.
Also, from the relation $K_Z(x,y)=Z(xy^*)=V(x)^{\sharp }V(y)$
and the linearity of $Z$, it follows that  
$\pi(a)\xi+\pi(b)\xi=\pi(a+b)\xi$
for all $\xi\in\cD(\pi)$ and $a,b\in\cA$.
Therefore, $\pi $ is a closable representation of $\cA$.
For $x,y,a\in \cA$ and $\lambda ,\mu \in \CC$, 
we have
$$\begin{array}{lll}
[\pi (a)V(x)\lambda ,V(y)\mu ]_{\cK}&=&
[V(xa^*)\lambda ,V(y)\mu ]_{\cK} \\[1em]
 &=&V(y)^{\sharp }V(xa^*)\lambda \overline{\mu } \\[1em]
 &=&K(y,xa^*)\lambda \overline{\mu } \\[1em]
 &=&K(ya,x)\lambda \overline{\mu } \\[1em]
 &=&[V(x)\lambda ,V(ya)\mu ]_{\cK} \\[1em]
 &=&[V(x)\lambda ,\pi (a^*)V(y)\mu ]_{\cK},
\end{array}$$
so that the domain of $\pi (a)^{\sharp } $  
contains $\cD (\pi )$ and \eqref{adj} holds.
Therefore, for every $a\in \cA$ the 
operator $\pi (a)$ is closable and $\pi $  
is a hermitian representation of the algebra $\cA$ on the Kre\u\i n
space $\cK$.
Also, for $a\in \cA$,
$$\begin{array}{lll}
Z(a)&=&K_Z(a,1)=V^{\sharp }(a)V(1) \\[1em]
&=&[V(1)1,V(a)1]_{\cK} \\[1em]
&=&[V(1)1,\pi (a^*)V(1)1]_{\cK} \\[1em]
&=&[\pi (a)V(1)1,V(1)1]_{\cK}.
\end{array}$$
Clearly, $\bigvee_{a\in\cA}\pi(a)V(1)1
=\bigvee_{x\in\cA}V(x){\CC}$, therefore, $(\pi ,\cK ,V(1)1)$ is a GNS data
associated to $Z$.

It is easy  now to see that the association of the GNS data
$(\pi,\cK,V(1)1)$, as in \eqref{pia}, to an arbitrary Kolmogorov
decomposition $(V,\cK)$ is a two-sided inverse to the correspondence
defined as in \eqref{kolGNS}.
\end{proof}

Proposition \ref{GNSkol} reduces the characterization of those hermitian 
functionals that admit GNS data to Theorem~\ref{kolmo}. A different
characterization was obtained in Theorem~2 in \cite{hof}.

\begin{theorem}\label{starex}
Let $\cA$ be a unital $*$-algebra and let $Z$ be a linear
hermitian
functional on $\cA$ with $Z(1)=1$.
Then $Z$ admits GNS data if and only if there exists a positive
definite scalar kernel $L$ on $\cA$
such that
\begin{equation}\label{was}
|Z(\sum_{i,j=1}^n\lambda_i\overline{\lambda}_j x_ix_j^*)|\leq 
\sum_{i,j=1}^n \lambda_i\overline{\lambda}_jL(x_i,x_j),\quad n\in\NN,\
\{\lambda_i\}_{i=1}^n\subset\CC,\ \{x_i\}_{i=1}^n\subset\cA . \end{equation}
\end{theorem}
\begin{proof} Note that \eqref{was} is equivalent to $-L\leq K_Z\leq
L$ and then apply Proposition~\ref{GNSkol} and Theorem~\ref{kolmo}.
\end{proof}

We now discuss the uniqueness property of the GNS data, an issue
previously addressed in \cite{hof}, but not
completely solved. Two GNS data $(\pi _1,\cK_1,\Omega_1)$ and
$(\pi _2,\cK _2,\Omega _2)$
are {\it unitarily equivalent} if there exists a unitary operator
$\Phi \in \cL (\cK _1, \cK _2)$ such that 
$\Phi \cD (\pi _1)=\cD(\pi _2)$, 
$\pi _2(a)\Phi =\Phi \pi _1(a)$
for all $a\in \cA$, and $\Phi \Omega _1=\Omega _2$.

\begin{theorem}\label{starun}
Let $\cA$ be a unital $*$-algebra and let $Z$ be a linear hermitian
functional on $\cA$ with $Z(1)=1$, admitting GNS data.
The following assertions are equivalent:

{\rm (1)} All GNS data of $Z$ are unitarily equivalent.

{\rm (2)} For each positive definite kernel $L$ on $\cA$
such that $-L\leq K_Z\leq L$, there exists $\epsilon >0$
such that either $(0,\epsilon )\subset \rho (A_L)$
or $(-\epsilon ,0)\subset \rho (A_L)$, where 
$A_L$ is the Gram operator of $K_Z$ with respect
to $L$.
\end{theorem}
\begin{proof} 
Let $(V_i,\cK _i)$, $i=1,2$, be two Kolmogorov decompositions
of $K_Z$ that are unitarily equivalent, that is, there exists
a unitary operator $\Phi \in \cL (\cK _1,\cK _2)$
such that $V_2(x)=\Phi V_1(x)$. Let $(\pi _i,\cK _i,\Omega _i)$,
$i=1,2$, be the corresponding GNS data for $Z$ as in 
Proposition~\ref{GNSkol}.
Then,
$$\cD (\pi _2)=\bigvee_{x\in\cA}V_2(x){\CC}=
\bigvee_{x\in\cA}\Phi V_1(x){\CC}=
\Phi (\bigvee_{x\in\cA}V_1(x){\CC})=\Phi \cD (\pi _1).$$
Also, for $a\in cA$ and $\lambda \in \CC$,
$$\pi _2(a)\Phi V_1(x)\lambda =
\pi _2(a)V_2(x)\lambda =V_2(xa^*)\lambda =
\Phi V_1(xa^*)\lambda =\Phi \pi _1(a)V_1(x)\lambda ,$$
which implies that $\pi _2(a)\Phi =\Phi \pi _1(a)$. Finally,
$$\Phi \Omega _1=\Phi V_1(1)1=V_2(1)1=\Omega _2,$$
therefore $(\pi _1,\cK _1,\Omega _1)$ and 
$(\pi _2,\cK _2,\Omega _2)$
are unitarily equivalent GNS data for $Z$.

Conversely, let $(\pi _i,\cK _i,\Omega _i)$,
$i=1,2$, be two unitarily equivalent GNS data for $Z$
and let $(V_i,\cK _i)$, $i=1,2$, be the Kolmogorov decompositions
of $K_Z$ associated to these GNS data by 
Proposition~\ref{GNSkol}. Therefore, there exists a unitary 
operator $\Phi \in \cL (\cK _1, \cK _2)$ such that 
$\Phi \cD (\pi _1)=\cD(\pi _2)$, 
$\pi _2(a)\Phi =\Phi \pi _1(a)$
for all $a\in \cA$ and $\Phi \Omega _1=\Omega _2$.
It follows that 
$$V_2(x)\lambda =\pi _2(a^*)\lambda \Omega _2=
\pi _2(a^*)\lambda \Phi \Omega _1=
\pi _2(a^*)\Phi \lambda \Omega _1=
\Phi \pi _1(a^*)\lambda \Omega _1=\Phi V_1(x)\lambda ,$$
which shows that $(V_1,\cK _1)$ and $(V_2,\cK _2)$
are unitarily equivalent Kolmogorov decompositions of the 
kernel $K_Z$. Now, an application of Theorem~\ref{unic}
concludes the proof.
\end{proof}

Another consequence of the Kolmogorov decomposition approach 
is the possibility of obtaining a characterization of those hermitian
functionals 
$Z$ that admit {\it bounded GNS data}, that is, the representation 
$\pi$ is made of bounded operators.

\begin{theorem}\label{starmar}
Let $\cA$ be a unital $*$-algebra and let $Z$ be a linear
hermitian functional on $\cA$ with $Z(1)=1$. Then $Z$ admits bounded
GNS data if and only if there exists a positive definite scalar kernel $L$ on
$\cA$ having the property \eqref{was} and such that 
for every $a\in\cA$ there exists $C_a>0$
with the property that
$$\sum_{i,j=1}^n \lambda_i\overline{\lambda}_jL(x_ia^*,x_ja^*)\leq C_a 
\sum_{i,j=1}^n \lambda_i\overline{\lambda}_jL(x_i,x_j),\quad  n\in\NN,\
\{\lambda_i\}_{i=1}^n\subset\CC,\ \{x_i\}_{i=1}^n\subset\cA.$$

\end{theorem}
\begin{proof} This is a consequence of Theorem~\ref{baza} and 
Proposition~\ref{GNSkol}.
\end{proof}

We conclude this section with a discussion of the
Jordan decomposition of a linear hermitian functional
on $\cA$, that is, the possibility of writing the hermitian functional
as the difference of two positive functionals. 
Let us first note that a functional $F:\cA\rightarrow\CC$ is
{\it positive}, that is, $F(a^*a)\geq 0$ for all $a\in\cA$, if and only if
the kernel $K_F$ associated to $F$ by the formula
\eqref{nucstar}
is positive definite. Also, if $F$ is a positive functional on $\cA$, then
$K_F$ is $\phi$-bounded, with the action $\phi$ defined as in
\eqref{actstar},
if and only if for any $a\in\cA$ there exists $C_a>0$ such that
\begin{equation}\label{fmarg} F(xa^*ax^*)\leq C_a F(xx^*),
\quad x\in\cA.\end{equation}
For simplicity, we call the positive functional $F$
$\phi$-{\it bounded} if $K_F$ is $\phi$-bounded.

Let $F_1$, $F_2$ be two positive functionals on the $*$-algebra
$\cA$. Then $F_1\leq F_2$, by definition, if $F_2-F_1$ is a positive
functional. It is easy to see that $F_1\leq F_2$ if and only if $K_{F_1}\leq
K_{F_2}$. The functionals $F_1$ and $F_2$
are called {\it disjoint} if their associated kernels
$K_{F_1}$ and $K_{F_2}$ are disjoint.

\begin{theorem}\label{starmar2}
Let $\cA$ be a unital $*$-algebra, let $Z$ be a linear
hermitian functional on $\cA$ with $Z(1)=1$, and let
$\phi $ be the action given by \eqref{actstar}.
The following assertions are equivalent:

{\rm (1)} There exists a linear positive $\phi$-bounded 
functional $Z_0$ on $\cA$
such that $-Z_0\leq Z\leq Z_0$.

{\rm (2)} $Z$ admits bounded GNS data $(\pi ,\cK ,\Omega )$ such that
the representation $\pi $ is similar with a $*$-representation on a
Hilbert space.

{\rm (3)} $Z=Z_+-Z_-$ for two disjoint linear positive definite functionals
on $\cA$ with the property that $(Z_++Z_-)$ is $\phi$-bounded.
\end{theorem}
\begin{proof}
The implications $(1)\Rightarrow (2)\Rightarrow (3)$
are direct consequences of Theorem~\ref{balabusta}
and Proposition~\ref{GNSkol}. For $(3)\Rightarrow (1)$
we use the proof of Theorem~\ref{baza}
in order to deduce that there
exists $L\in \cB _{\phi }^+(\cA,\CC)$ such that 
$-L\leq K_Z\leq L$. Then  
Theorem~\ref{balabusta} shows that 
$L(x,\phi (a,y))=L(\phi (a^*,x),y)$
for all $x,y,a\in \cA$. Also, in this case, 
$L$ is linear in the first variable (hence, antilinear in the 
second variable). If we define 
$$Z_0(x)=L(x,1)$$
for $x\in \cA$, then $Z_0$ is a linear functional on $\cA$ and 
$$K_{Z_0}(x,y)=Z_0(xy^*)=L(xy^*,1)=L(x,y).$$
Now all the required properties of $Z_0$ follow from the 
corresponding properties of $L$.
\end{proof}

\begin{remark} {\rm It is interesting to note that under fairly general
assumptions on the $*$-algebra $\cA$, every positive functional $F$ on
$\cA$ is $\phi$-bounded, that is, for all $a$ in $\cA$ we have
\eqref{fmarg}. This holds, for instance, if $\cA$ is a Banach
$*$-algebra, cf. Lemma 37.6 in \cite{BD}, with the constant $C_a$ equal
to the spectral radius of $a^*a$.}\end{remark}

\section{Dilation theory}

We can extend the construction of bounded GNS data
to the case of hermitian mappings, a topic usually
referred to as dilation theory. In this section
we obtain a version of the Stinespring theorem and 
we discuss the connection between invariant Kolmogorov decompositions and
the Wittstock's result on representing hermitian completely bounded 
maps as the difference of two completely positive maps.

Let $\cA$ be a unital $*$-algebra and let $T:\cA\rightarrow \lh$ be a linear
hermitian map, where $\cH$ is a Kre\u\i n space. A {\it Stinespring
dilation} of $T$ is, by definition, a triple $(\pi,\cK,B)$ where,
$\cK$ is a Kre\u\i n space, $\pi\colon\cA\rightarrow \cL(\cK)$ is a
selfdajoint representation, and $B\in\cL(\cH,\cK)$ such that
\begin{equation}
T(a)=B^{\sharp }\pi (a)B,\quad a\in\cA.\end{equation}
If, in addition, $\bigvee_{a\in\cA}\pi(a)B\cH$ is dense in $\cK$, then
the Stinespring dilation is called {\it minimal}.

We consider the action $\phi$ of $\cA$ on itself defined
as in the previous section by the formula \eqref{actstar},
and a hermitian kernel is associated to $T$ by the formula
\begin{equation}K_T(x,y)=T(xy^*),\quad x,y\in\cA.\end{equation}
It follows readily that $K_T$ is $\phi$-invariant, that is,
\begin{equation}\label{siminica}
K_T(x,\phi (a,y))=K_T(\phi (a^*,x),y).
\end{equation}

\begin{lemma}\label{stinekol} Given a minimal Stinespring dilation 
$(\pi,\cK,B)$ of the hermitian linear map $T\colon\cA\rightarrow\cL(\cH)$,
let
\begin{equation}\label{vepi}V(x)=\pi(x^*)B,\quad x\in\cA.\end{equation}
Then $(V,\cK)$ is an invariant Kolmogorov decomposition of the
hermitian kernel $K_T$. In addition, \eqref{vepi} establishes a
bijective correspondence between the set of minimal Stinespring
dilations of $T$ and the set of invariant Kolmogorov decompositions of
$K_T$.\end{lemma} 
\begin{proof} Let $(\pi,\cK,B)$ be a minimal Stinespring dilation of
$T$ and define $(V,\cK)$ as in \eqref{vepi}. Then
$$V(x)^\sharp V(y)=B^\sharp \pi(x^*)^\sharp \pi(y^*)B=B^\sharp
\pi(xy^*)B=T(xy^*)=K_T(x,y),\quad x,y\in\cA.$$
Since $\bigvee_{a\in\cA}\pi(a)B\cH=\bigvee_{x\in\cA}V(x)\cH$ it
follows that $(V,\cK)$ is a Kolmogorov decomposition of $T$. Let us
note that, by the defintion of $V$,
$$\pi(a)V(x)=\pi(a)\pi(x^*)B=\pi(ax^*)B=V(xa^*)\,\quad a,x\in\cA,$$
and hence, letting $U=\pi$, it follows that the Kolmogorov
decomposition $(V,\cK)$ is invariant.

Conversely, let $(V,\cK)$ be an invariant Kolmogorov decomposition of the
hermitian kernel $K_T$, that is, there exists a hermitian
representation $U\colon\cA\rightarrow\cL(\cK)$ 
of the multiplicative semigroup with involution $\cA$, such that
$$U(a)V(x)=V(xa^*),\quad a,x\in\cA.$$
Define $\pi=U$ and $B=V(1)$. Since $T$ is linear it follows easily
that $\pi$ is also linear, hence a selfadjoint representation of the
$*$-algebra $\cA$ on the Kre\u\i n space $\cK$. Then, taking into account
that $V(a)=U(a^*)B$ for all $a\in\cA$, it follows
$$T(a)=V(a)^\sharp V(1)=B^\sharp U(a^*)^\sharp B=B^\sharp U(a)B,\quad
a\in\cA,$$ and since
$\bigvee_{a\in\cA}\pi(a)B\cH=\bigvee_{x\in\cA}V(x)\cH$ we thus proved
that $(\pi,\cK,B)$ is a minimal Stinespring dilation of $T$.

We leave to the reader to show that the mapping defined in
\eqref{vepi} is a two-sided inverse of the mapping associating to each
invariant Kolmogorov decomposition $(V,\cK)$ the minimal Stinespring
dilation $(\pi,\cK,B)$ as above.\end{proof}

As a consequence of Theorem~\ref{baza} and Lemma~\ref{stinekol} we obtain.

\begin{theorem}\label{stine}
Let $\cA$ be a unital $*$-algebra, $\cH$ a Kre\u\i n space, 
and let $T:\cA\rightarrow \lh$ 
be a linear hermitian map. The following assertions are equivalent:

{\rm (1)} There exists a positive definite kernel 
$L\in \cB _{\phi }^+(\cA,\cH)$, $\phi $ given by \eqref{actstar}, 
such that $-L\leq K_T\leq L$.

{\rm (2)} $T$ has a minimal Stinespring dilation. 
\end{theorem}

In case $\cA$ is a $C^*$-algebra and $\cH$ is a Hilbert space, more
can be said. Under these assumptions, recall 
first that for $n\in \NN$, $M_n(\cA)=M_n\otimes
\cA$ is a
$C^*$-algebra in a natural (and unique) way, where $M_n$ denotes the
$C^*$-algebra of matrices with complex entries of size $n\times n$. 
Recall also that a 
map $S\colon\cA \rightarrow \lh $ is called {\it completely positive} if
for all $n\in\NN$ and any positive matrix $[a_{ij}]_{i,j=1}^n\in
M_N(\cA)$, the matrix $[S(a_{ij})]_{i,j=1}^n$ is positive in
$M_n(\lh)$. Since every positive matrix $[a_{ij}]_{i,j=1}^n$ with
entries in $\cA$ can be written as a sum of matrices of the form
$[x_i^*x_j]_{i,j=1}^n$ with $x_i\in \cA$ for $i=1,\ldots,n$, it
follows that the map $S\colon\cA \rightarrow \lh $
is completely positive if and only if the kernel $K_S$ is positive definite.
Also, we say that two completely positive maps
$F_1$ and $F_2$ are {\it disjoint} if the
associated kernels $K_{F_1}$ and $K_{F_2}$ are disjoint.

The following result is the Wittstock's characterization of completely bounded
maps, cf.\ \cite{W}. Recall that, if $\cA$ is a $C^*$-algebra and
$\cH$ is a Hilbert space, a linear map $T\colon\cA\rightarrow\cL(\cH)$
is called {\it completely bounded} if there exists $C>0$ such
for all $n\in\NN$ and $[a_{ij}]_{i,j=1}^n\in M_n\otimes \cA$  we have
$\|[T(a_{ij})]_{i,j=1}^n\|\leq C\|[a_{ij}]^n_{i,j=1}\|$, where the
norms are calculated in the $C^*$-algebras $M_n\otimes\lh$ and,
respectively, $M_n\otimes\cA$.

\begin{theorem}\label{witt}
Let $\cA$ be a unital $C^*$-algebra and let $T:\cA \rightarrow \lh$
be a linear hermitian map. The following assertions are equivalent:

{\rm (1)} $T$ is completely bounded.

{\rm (2)} There exists a completely positive map $S:\cA \rightarrow \lh$
such that $-S\leq T\leq S$.

{\rm (3)} $T$ has a minimal
Stinespring dilation similar to a (minimal) Stinespring dilation on
a Hilbert space. 
 
{\rm (4)} $T=T_+-T_-$ for two disjoint completely positive maps $T_+$ and 
$T_-$.
\end{theorem}
\begin{proof}
In the following we let $\cU(\cA)$ be the unitary group of $\cA$.
Then $\cU(\cA)$ has the involution $\fI (a)=a^{-1}=a^*$
and acts on $\cA$ by $\phi (a,x)=xa^*=xa^{-1}$.

$(1)\Rightarrow (2)$ 
 We use the
Paulsen's off-diagonal technique. Briefly, 
assume that $T$ is completely bounded. By Theorem~7.3 in 
\cite{paul}, there exist completely positive maps
$\phi _1$ and $\phi _2$ such that the map
$$F(\left[\begin{array}{cc} a & b\\ c &
d
\end{array}\right])=
\left[\begin{array}{cc}
 \phi _1(a) & T(b) \\
 T(c^*)^* & \phi _2(d) 
\end{array}\right]$$
is completely positive. Define $S(a)=1/2(\phi _1(a)+\phi _2(a))$, 
which is a completely positive map. We can check that 
$-S\leq T\leq S$. First, let $a\geq 0$, $a\in \cA$.
Then 
$\left[\begin{array}{cc}
a & \pm a \\
 \pm a & a 
\end{array}\right]\geq 0$,
so that
$\left[\begin{array}{cc}
\phi _1(a) & \pm T(a)\\
\pm T(a) & \phi _2(a)
\end{array}\right]\geq 0$. In particular, for $\xi \in \cH$,
$$\langle \left[\begin{array}{cc}
\phi _1(a) & \pm T(a)\\
\pm T(a) & \phi _2(a)
\end{array}\right]
\left[\begin{array}{cc} 
\xi \\ \xi \end{array}\right],
\left[\begin{array}{cc} 
\xi \\ \xi \end{array}\right]
\rangle \geq 0,$$
or $\langle (\phi _1(a)\pm 2T(a)+\phi _2(a))\xi ,\xi\rangle \geq 0$.
Therefore, $S\pm T$ are positive maps. The argument can
be extended in a straightforward manner (using the so-called
canonical shuffle as in \cite{paul}) to show that 
$S\pm T$ are completely positive maps.

$(2)\Rightarrow (3)$ Since $S$ is completely positive, the 
kernel $K_S$ is positive definite and satisfies 
$-K_S\leq K_T\leq K_S$. Also, 
$$K_S(x,\phi (a,y))=K_S(\phi (a^{-1},x),y),\quad a\in\cU(\cA),\ x,y\in\cA.$$ 
By Theorem~\ref{corolar}, there exists a 
Kolmogorov decomposition $(V,\cK )$
of $K_T$ and a fundamentally reducible representation
$U$ of $\cU(\cA )$ on $\cK $ such that
$$U(a)V(x)=V(xa^{-1}),\quad a\in\cU(\cA),\ x\in\cA.$$
 Let $J$ be a fundamental 
symmetry on $\cK$ such that $U(a)J=JU(a)$ for all $a\in \cU(\cA)$.
Then $U$ is also a representation of $\cU(\cA )$ on the Hilbert
space $(\cK ,\langle \cdot ,\cdot \rangle _J)$. Also, for $a\in \cU(\cA)$,
$$T(a)=K_T(a,1)=V(a)^{\sharp }V(1)=V(1)^*JU(a)V(1).$$
Since $\cA$ is the linear span of $\cU(\cA )$ and $T$ is linear, 
$U$ can be extended by linearity to a representation $\pi $ 
of $\cA$ on $\cK$ commuting with $J$ and such that  
$$T(a)=V(1)^*J\pi (a)V(1)$$
holds for all $a\in \cA$. Also, $\bigvee _{a\in \cU(\cA )}U(a)V(1)\cH
=\bigvee _{a\in \cU(\cA )}V(a^{-1})\cH =\bigvee _{a\in \cU(\cA )}V(a)\cH $
and using once again the fact that $\cA$ is the linear span of $\cU(\cA)$, 
we deduce that $\bigvee _{a\in \cA }U(a)V(1)\cH $ is dense in $\cK$. 
 
(3)$\Rightarrow$(4) Note first that statement (3) is equivalent with
the following: there exists a Hilbert space $\cK$, a $*$-representation 
$\pi :\cA\rightarrow \lk $, and a bounded operator $B\in \lhk $
such that
$$T(a)=B^*J\pi (a)B,\quad a\in \cA,$$
where $J$ is a symmetry on $\cK$
commuting with $\pi (a)$ for all $a\in \cA$, and 
$\bigvee _{a\in \cA}\pi (a)B\cH$ is dense in $\cK$.

Assuming this, define for $a\in \cA $,
$$V(a)=\pi (a^*)B.$$
Then $V(a)$ is in $\lhk $ and 
$\bigvee _{a\in \cA}V(a)\cH =\bigvee _{a\in \cA}\pi (a)B\cH $.
$\cK $ becomes a Kre\u{\i}n space by setting $[x,y]_{\cK }=
\langle Jx,y\rangle $, $x,y\in \cK $.
Also, for $\xi ,\eta \in \cH $,
$$\begin{array}{lll}
\langle V(x)^{\sharp }V(y)\xi ,\eta \rangle &=&
\langle JV(y)\xi ,V(x)\eta \rangle \\[1em]
&=& \langle J\pi (y^*)B\xi ,\pi (x^*)B\eta \rangle \\[1em]
&=& \langle J\pi (xy^*)B\xi ,B\eta \rangle \\[1em]
&=& \langle T(xy^*)\xi ,\eta \rangle \\[1em]
&=& \langle K_T(x,y)\xi ,\eta \rangle ,
\end{array}$$
so that $(V,\cK )$ is a Kolmogorov decomposition of $K_T$.
Let $J=J_+-J_-$ be the Jordan decomposition
of $J$ and define, as in the proof of Theorem~\ref{baza},
the hermitian kernels
$$K_{\pm }(x,y)=\pm V(x)^{\sharp }J_{\pm }V(y).$$
It was checked in the proof of Theorem~\ref{corolar}
that 
$$K_{\pm }(x,\phi (a,y))=K_{\pm }(\phi (a^{-1},x),y)$$
for all $x,y\in \cA $ and $a\in \cU(\cA )$. For $x\in \cA $ define
$$T_{\pm }(x)=K_{\pm }(x,1).$$
Then $T_{\pm }(x)=\pm B^{\sharp }\pi (x^*)^{\sharp }J_{\pm }V(1)$
are linear maps on $\cA$ and for $x\in \cA$, $y\in \cU(\cA )$, 
we get
$$\begin{array}{lll}
K_{T_{\pm }}(x,y)&=&T_{\pm }(xy^{-1}) \\[1em]
&=&K_{\pm }(xy^{-1},1) \\[1em]
&=&K_{\pm }(\phi (y,x),1) \\[1em]
&=&K_{\pm }(x,\phi (y^{-1},1)) \\[1em]
&=&K_{\pm }(x,y).
\end{array}$$
Since $K_{T_{\pm }}$ and $K_{\pm }$ are antilinear in the
second variable and $\cA $ is the linear span of $\cU(\cA )$,
it follows that
$$K_{T_{\pm }}(x,y)=K_{\pm }(x,y)$$
for all $x,y\in \cA$. This implies that
$T_{\pm }$ are disjoint completely positive maps 
such that $T=T_+-T_-$.

The implication (4)$\Rightarrow$(1) follows from Theorem~\ref{corolar}.
\end{proof}

\begin{remark} {\em We can also notice that 
it follows easily from the proof of Theorem~\ref{witt}
that the Wittstock representation in the statement (3) is 
unique up to unitary equivalence for completely bounded 
hermitian maps. 
Since the representation of a completely bounded 
map as a sum of two hermitian completely bounded maps is not 
unique, the uniqueness property is lost in the general nonhermitian
case (see \cite{paul}).
}
\end{remark}

\end{document}